\documentclass[oneside,english,reqno]{amsart}
\usepackage{geometry}

\usepackage{amsthm}
\usepackage{amssymb}
\usepackage{setspace}
\usepackage{tikz}
\usepackage{hyperref}
\usepackage{enumitem}
\usepackage{bbm}
\usetikzlibrary{arrows}
\tikzstyle{vertex}=[draw, circle, scale=0.5]
\tikzstyle{edge}=[draw=white, circle, scale=0.7, fill=white]
\usepackage{caption}
\usepackage{subcaption}
\usetikzlibrary{calc}

\makeatletter
\numberwithin{equation}{section}
\numberwithin{figure}{section}
\theoremstyle{plain}
\newtheorem{theorem}{Theorem}
\newtheorem*{theorem*}{Theorem}
  \newtheorem{cor}[theorem]{Corollary}
  \newtheorem{lemma}[theorem]{Lemma}
  \newtheorem{defn}[theorem]{Definition}
  \newtheorem{conj}[theorem]{Conjecture}
  \newtheorem*{conj*}{Conjecture}
  \newtheorem{prop}[theorem]{Proposition}
  \newtheorem{obs}[theorem]{Observation}
  \newtheorem{question}[theorem]{Question}

  \newtheorem{remark}[theorem]{Remark}

\numberwithin{theorem}{section}

\newcommand{\C}{\mathcal{C}}
\newcommand{\Z}{\mathbb{Z}}

\newcommand{\E}{\operatorname{\mathbb{E}}}

\newcommand{\rk}{\mathcal{E}}
\newcommand{\fh}{\mathcal{H}}
\newcommand{\fH}{\mathcal{H}}
\newcommand{\ex}{\operatorname{ex}}

\renewcommand{\O}{\mathcal{O}}
\renewcommand{\backslash}{\setminus}
\renewcommand{\deg}{d}
\renewcommand{\subset}{\subseteq}
\newcommand{\cal}[1]{\mathcal{#1}}

\begin{document}

\title{Inverting the Tur\'an problem}

\author{Joseph Briggs$^1$ \and Christopher Cox$^1$}

\begin{abstract}
Classical questions in extremal graph theory concern the asymptotics of $\ex(G, \fH)$ where $\fh$ is a fixed family of graphs and $G=G_n$ is taken from a ``standard'' increasing sequence of host graphs $(G_1, G_2, \dots)$, most often $K_n$ or $K_{n,n}$. Inverting the question, we can instead ask how large $e(G)$ can be with respect to $\ex(G,\fH)$. We show that the standard sequences indeed maximize $e(G)$ for some choices of $\fh$, but not for others. Many interesting questions and previous results arise very naturally in this context, which also, unusually, gives rise to sensible extremal questions concerning multigraphs and non-uniform hypergraphs.
\end{abstract}

\maketitle

\footnotetext[1]{Carnegie Mellon University, Pittsburgh, PA, USA. \texttt{\{jbriggs,cocox\}@andrew.cmu.edu}.}

\section{Introduction and Motivation}
For a graph $G$ and a family of graphs $\fh$, the extremal number of $\fh$ in $G$ is defined to be
\[
\ex(G, \fh) = \max \bigl\{e(F):  F \subseteq G\text{ and }  H \not\subseteq F \text{ for any } H \in \fh\bigr\}.
\]
When the family consists only of a single graph, $\ex(G,H)$ is used in place of $\ex(G,\{H\})$.

A typical example of this is when $\fh =\{C_3, C_4, C_5, \dots \}$ is the collection of all cycles, in which case the extremal number is simply the \emph{graphic matroid rank} of $G$, an important graph parameter in its own right.

The Tur\'{a}n problem, one of the cornerstones of extremal graph theory concerns the behavior of $\ex(K_n,H)$ for a fixed $H$ when $n$ is large. The first result along these lines is a theorem of Mantel (see, for instance~\cite{Bollobas13}) which states that $\ex(K_n,K_3)=\lfloor n^2/4\rfloor$. Tur\'an~\cite{Turan41} obtained a version for $K_t$ in place of $K_3$, in particular obtaining $\ex(K_n,K_t)=\bigl(1-{1\over t-1}+o(1)\bigr){n^2\over 2}$ where $o(1)\to 0$ as $n\to\infty$. In a similar spirit, the Erd\H{o}s--Stone Theorem~\cite{ES46} states that if $\chi=\chi(H)$ is the chromatic number of $H$, then $\ex(K_n,H)=\bigl(1-{1\over\chi-1}+o(1)\bigr){n^2\over 2}$. The Erd\H{o}s--Stone Theorem asymptotically answers the Tur\'an problem, except when $H$ is bipartite, in which case the bound becomes $o(n^2)$. In this situation, known as the degenerate case, the asymptotic behavior of very few graphs is known and is an active area of research (c.f.~\cite{FS13}).

Most approaches in the case of a bipartite graph instead ask about $\ex(K_{n,n},H)$, which is known as the Zarankiewicz problem~\cite{Zaran51}. This is often seen as a more natural question and provides bounds on the Tur\'an problem as ${1\over 2}\ex(K_n,H)\leq\ex(K_{n/2,n/2},H)\leq\ex(K_n,H)$ for bipartite $H$. In the special case of $H=C_4$, the incidence graphs showing tightness for the Zarankiewicz problem were spotted a few years before polarity graphs showing tightness for the Tur\'an problem (see
~{\cite[Section 3]{FS13}}).

With this in mind, we set out to explore a framework in which to ask: what is the most ``natural'' or ``best'' host graph for a fixed family of graphs? This suggests optimizing a particular monotone graph parameter over all host graphs $G$ where $\ex(G,\fH)$ is bounded, the simplest of which is just the edge count. Thus we define the following extremal function for $\fh$:
\[
\rk_\fh(k):=\sup\bigl\{e(G):\ex(G,\fh) < k\bigr\}.
\]
In other words, for a family $\fh$, we would like to determine the host graph $G$ with the most edges such that any $k$ edges from $G$ \emph{must} contain a copy of some $H\in\fh$. In other words, $G$ is best at ``forcing'' a copy of some $H\in\fh$.
When the family consists only of a single graph, we write $\rk_H(k)$ in place of $\rk_{\{H\}}(k)$. Note that it is necessary to consider the supremum here as $\rk_{\fh}(k)$ may be infinite. In particular, $\rk_{K_{1,t}}(k)=\rk_{tK_2}(k)=\infty$ for $k\geq t$ as for any $s\geq t$, $\ex(K_{1,s},K_{1,t})=t-1=\ex(sK_2, tK_2)$, despite both host graphs having $s$ edges. However, we will show in Proposition~\ref{prop:sunflower} that stars and matchings classify all families having $\rk_\fh(k)=\infty$.

In a similar fashion to the original Tur\'an problem, this paper considers two questions:
\begin{itemize}
\item What are the asymptotics of $\rk_\fh(k)$?
\item When $\rk_\fh(k)$ can be determined precisely, which host graphs $G$ attain $e(G)=\rk_\fh(k)$?
\end{itemize}

On the one hand, we will show that for nonbipartite $H$, this question behaves more or less as one might expect.  For example, the following theorem is close in spirit to the Erd\H{o}s--Stone Theorem:
\begin{theorem}\label{thm:esrank}
If $\fH$ is a family of graphs with $\rho=\min\{\chi(H): H \in \fH\}\geq 3$, then 
\[
\rk_\fH(k)=\biggl(1+{1\over\rho-2}+o(1)\biggr)k,
\]
where $o(1)\to 0$ as $k\to\infty$.
\end{theorem}

Recalling our motivation from the Zarankiewicz problem, we show that complete bipartite graphs are optimal hosts for at least one natural family, namely the collection $\C_e:=\{C_4,C_6,\dots\}$ of even cycles:
\begin{theorem}\label{thm:evencycles}
For $k \geq 4$,
$\rk_{\C_e}(k) = \bigl\lfloor \frac{k^2}{4} \bigr\rfloor$, with $K_{\lfloor k/2\rfloor,\lceil k/2\rceil}$ being the unique extremal graph for $k\geq 6$.
\end{theorem}

On the other hand, this is already a challenge for the case $H=K_{2,2}$:
\begin{question}
What is $\rk_{C_4}(k)$ and what is the optimal host graph?
\end{question}

One peculiar feature of our question is that it is sensible even for multigraphs or nonuniform hypergraphs.  We let $\rk^*_\fh(k)$ denote the maximum number of edges among host \emph{multigraphs} $G$ with $\ex(G,\fh) < k$. Here a multigraph is a hypergraph with parallel edges. The parameter $\rk^*_\fh(k)$ will be important in proving bounds on $\rk_\fh(k)$ when $\fh$ is a family of simple graphs. However, we do not even know the following:

\begin{conj}\label{conj:multi}
If $\fh$ is a family of simple $r$-uniform hypergraphs, then $\rk_\fH(k) = \rk^*_\fh(k)$.
\end{conj}

Curiously, for \emph{non-uniform} graphs $H$ without parallel edges, that is $H$ which have edges of different uniformities yet each edge still has multiplicity 1, the analogue of the above conjecture fails:

\begin{theorem}\label{thm:dumbbell}
Let $\cal{D}$ denote the non-uniform hypergraph with a single edge and a loop at each vertex (the dumbbell graph). Then $\rk_{\cal{D}}(k) = \bigl\lfloor\frac{3}{2}(k-1)\bigr\rfloor$, whereas $\rk^*_{\cal{D}}(k) \sim \phi\cdot k$ where $\phi=1.618\ldots$ is the golden ratio.
\end{theorem}

Beyond the theorems stated above, in our study of $\rk_\fH(k)$ and optimal host graphs, we will also show that:
\begin{enumerate}
\item \label{itm:unifinf} For uniform hypergraphs $H$, $\rk_H(k)$ is only infinite for sunflowers (Proposition~\ref{prop:sunflower}),
\item \label{itm:cliques} Cliques are best at forcing cliques (Theorem~\ref{thm:cliques}),
\item \label{itm:cycle} Cliques are best at forcing a cycle (Theorem~\ref{thm:cycles}),
\item \label{itm:p3k3} Complements of matchings are best at forcing $\{P_3,K_3\}$ (Theorem~\ref{thm:p3k3}),
\item \label{itm:p3} Cliques with pendant edges are best at forcing $P_3$ (Theorem~\ref{thm:p3}),
\item \label{itm:p1up2}Two disjoint cliques or a modified power of a cycle, depending on parity, are best at forcing $P_1 \cup P_2$ (Corollary~\ref{cor:p2p3} \& Theorem~\ref{thm:degmatch}),
\item \label{itm:1unif} For 1-uniform multigraphs $H$, $\rk_H^*(k)$ is quadratic in $k$ (Theorem~\ref{thm:1u}).
\end{enumerate}
In fact, for items \ref{itm:cliques}, \ref{itm:p3k3}, and \ref{itm:p1up2}, the correct  behavior of $\rk_{\fh}(k)$ is implicit in references \cite{alon1996bipartite}, \cite{FHHM02} and \cite{AHS72}, respectively, but our results will prove uniqueness of the respective host graphs.

This paper begins by exploring basic properties of $\rk_\fh(k)$ in Section~\ref{sec:basic}

The organization of this paper is as follows. In Section~\ref{sec:basic}, we begin our study of $\rk_\fh(k)$ and $\rk_\fh^*(k)$ by obtaining the natural analogue of the Erd\H{o}s--Stone theorem and then discussing the relationship between $\rk_\fh(k)$ and $\rk_\fh^*(k)$ in Section~\ref{sec:multi}. We then, in Section~\ref{sec:graphs}, determine $\rk_{K_t}(k),\rk_{K_t}^*(k)$ exactly for infinitely many values of $k$, and also look at $\rk_\fh(k)$ when $\fh$ is a family of cycles and when $\fh$ consists of small graphs, in some cases extending to $\rk_\fh^*(k)$. In Section~\ref{sec:nonu}, we look at $\rk_\fh(k),\rk_\fh^*(k)$ when the graphs in $\fh$ are not required to be uniform, and finally conclude by determining the behavior of $\rk_H^*(k)$ when $H$ is a 1-uniform multigraph in Section~\ref{sec:1u}. A list of conjectures and future directions can be found in Section~\ref{sec:conclusion}.

\subsection{Notation}

We follow standard notation from~\cite{West96} with the exception that $P_t$ will denote the path of $t$ \emph{edges}. For a (multi)(hyper)graph $G$, we let $V(G)$ denote the vertex set of $G$, $E(G)$ denote the edge set, and use $e(G)= |E(G)|$. Furthermore, throughout the paper, unless explicitly stated otherwise, all graphs will be assumed to have no isolated vertices. We will also abuse our language slightly and often refer to both multigraphs and hypergraphs simply as ``graphs.'' We will ensure that the type of graph discussed is always clear from context. Additionally, throughout the majority of the paper, all graphs will be assumed to be uniform, the only exception being in Section~\ref{sec:nonu}.

For a graph $G=(V,E)$ and disjoint subsets $S,T\subseteq V$, we use $G[S]$ to denote the subgraph of $G$ induced by $S$ and $G[S,T]$ to denote the subgraph of $G$ with vertex set $S\cup T$ where $xy\in E(G[S,T])$ if and only if $xy\in E$ and $x\in S$ and $y\in T$. In these cases we will use $e[S]=e(G[S])$ and $e[S,T]=e(G[S,T])$ for brevity.

For a graph $G$ and integer $t$, we denote the graph consisting of $t$ vertex-disjoint copies of $G$ by $tG$, e.g.\ $tK_2$ is the matching on $t$ edges. For integers $m\leq n$, we use $[m,n]=\{m,m+1,\dots,n\}$ and $[n]=[1,n]$.

\section{Basic Results}\label{sec:basic}

Recall that for a family of simple (hyper)graphs $\fh$, we defined
\[
\rk_\fh(k):=\sup\bigl\{e(G):\ex(G,\fh)<k\bigl\}.
\]
Also, for a family of (multi)(hyper)graphs $\fh$, we defined
\[
\rk_\fh^*(k):=\sup\bigl\{e(G):\ex(G,\fh)<k,\text{ $G$ is a multigraph}\bigr\}.
\]
Certainly if $\fh$ is a family of simple graphs, then we always have $\rk_\fh(k)\leq\rk_\fh^*(k)$ and we will sometimes rely on $\rk_\fh^*(k)$ to prove upper bounds on $\rk_\fh(k)$. This will be key in the proof of Theorem~\ref{thm:cliques}. We believe that it should be the case that $\rk_\fh(k)=\rk_\fh^*(k)$ if $\fh$ is a family of simple graphs, and we will discuss this question in Section~\ref{sec:multi}. 

Before addressing any questions about $\rk_\fh(k)$ and $\rk_\fh^*(k)$, we begin with a simple observation about two monotonicity properties which will be used extensively throughout the paper.

\begin{obs}
	If $H$ is a subgraph of $G$, then $\rk_H^*(k)\geq\rk_G^*(k)$ and if $\fh,\cal{G}$ are families of graphs with $\fh\subseteq\cal{G}$, then $\rk_\fh^*(k)\leq\rk_{\cal{G}}^*(k)$ (note the reversed inequality). Furthermore, if all graphs are simple, these same inequalities hold for $\rk_\fh(k)$.
\end{obs}

We begin our study of $\rk_\fh(k)$ by classifying exactly when $\rk_\fh(k)=\infty$. As we saw in the introduction, for 2-uniform graphs, this happens for stars and matchings. When considering 2-uniform multigraphs, we obtain an additional example: an edge with multiplicity $t$. Indeed, these examples classify precisely when $\rk_\fh(k),\rk_\fh^*(k)=\infty$ for 2-uniform (multi)graphs. For families of hypergraphs of higher uniformity, this question is answered by the classical sunflower lemma due to Erd\H{o}s and Rado~\cite{ER60}.

\begin{defn} Let $H$ be an $r$-uniform (multi)graph. $H$ is said to be a sunflower if for some $S\subseteq V(H)$, called the \emph{core} of $H$, every pair of distinct edges $e_1,e_2\in E(H)$ has $e_1\cap e_2=S$.
\end{defn}

Note that $K_{1,t}$ and $tK_2$ are precisely the simple 2-uniform sunflowers while the multiedge is a 2-uniform (multi)sunflower.

\begin{prop}\label{prop:sunflower}
	If $\fh$ is a family of $r$-uniform graphs, then $\rk_\fh^*(k)=\infty$ for $k$ sufficiently large if and only if $\fh$ contains a sunflower. Furthermore, if $\fh$ is a family of simple $r$-uniform graphs, this result also holds for $\rk_\fh(k)$.
\end{prop}

\begin{proof}
	We give a proof only for $\rk_\fh^*(k)$ since the case of $\rk_\fh(k)$ is very similar. First suppose that $H\in\fh$ is a sunflower with $e(H)=m$ and core $S\subseteq V(H)$. Then any $r$-uniform sunflower $G$ with $s\geq m$ edges and core of size $|S|$ has $\ex(G,H)=m-1$. Hence, if $k\geq m$, we have $\rk_\fh^*(k)\geq\rk_H^*(k)\geq s$, so $\rk_\fh^*(k)=\infty$.
	
	On the other hand, suppose that $\fh$ does not contain a sunflower. By the Erd\H{o}s--Rado sunflower lemma~\cite{ER60}, if $G$ is an $r$-uniform (multi)graph with $e(G)>r!(k-1)^{r+1}$ \footnote{When $G$ is simple, this can be lowered to $r!(k-1)^r$. Furthermore, better bounds are known in general, but they are unnecessary here.}, then $G$ contains a sunflower $F$ with at least $k$ edges. Since $\fh$ has no sunflower, $F$ contains no hypergraph in $\fh$, so $\ex(G,\fh)\geq k$. In other words, $\rk_\fh^*(k)\leq r!(k-1)^{r+1}<\infty$.
\end{proof}

Now that we have an understanding of exactly when $\rk_\fh(k)$ and $\rk_\fh^*(k)$ can be infinite, we turn to understanding the connection with the usual Tur\'an number $\ex(K_n,\fh)$ when $\fh$ is a family of simple graphs.

\begin{defn}
	For an arbitrary $r$-uniform hypergraph family $\fH$ 
	we denote by 
	\[
	\biggl( \pi_n(\fH) := \frac{\ex(K_n^{(r)},\fH)}{\binom{n}{r}}\biggr)_{n\geq 1}
	\]
	the sequence of Tur\'an densities and denote the limiting density $\pi(\fH):=\lim_{n\to\infty}\pi_n(\fH)$.
	$\fH$ is said to be \emph{degenerate} if $\pi(\fh)=0$.
\end{defn}

Note that $\bigl(\pi_n(\fH)\bigr)_{n\geq 1}$ is a decreasing sequence of densities for any $\fh$ by averaging over subgraphs, so the limit always exists.
Furthermore, there is a standard classification:

\begin{prop}
	An $r$-uniform graph $H$ is degenerate if and only if it is $r$-partite. That is to say, we may $r$-color $V(H)$ so that each $e\in E(H)$ has 1 vertex of each color, or equivalently, $H \subseteq K^{(r)}_{\underbrace{t,t,\dots,t}_r}$ for some $t$.
\end{prop}

Indeed, for $H$ nondegenerate, $\pi(H) \geq r!/r^r$ since the balanced $r$-partite hypergraph
$K^{(r)}_{n/r , \dots, n/r}
\not\supseteq H$, otherwise
$\ex(K_n^{(r)},H) =o(n^r)$ is true by an induction on $r$, as was observed by Erd\H{o}s~\cite{Erdos64}. In fact, these easily generalize to families of $r$-uniform graphs; namely $\pi(\fH)=0$ if and only if $\fH$ contains a degenerate graph. See \cite{Keevash11} for a survey on the hypergraph Tur\'an problem.

\begin{theorem}\label{thm:eshypergraphs}
	If $\fH$ is a family of simple $r$-uniform hypergraphs not containing a degenerate graph, then
	\[
	\rk_\fh(k),\rk_\fh^*(k)=\biggl({1\over\pi(\fH)}-o(1)\biggr)k.
	\]
\end{theorem}

\begin{proof}
	\textit{Lower bound.} For a positive integer $k$, let $n$ be the largest integer for which $k>\pi_n(\fH){n\choose r}$. Since $\pi_n(\fH)=\pi(\fH)+o(1)$ and ${n+1\choose r}-{n\choose r}=O(n^{r-1})$, we observe that $\pi_{n+1}(\fH){n+1\choose r}-\pi_n(\fH){n\choose r}= o(n^r)$; thus $k\leq\pi_n(\fH){n\choose r}+o(k)$. Then, as $\ex(K_n^{(r)},\fH)=\pi_n(\fH){n\choose r}<k$,
	\[
	\rk_\fH(k)\geq e(K_n^{(r)})={n\choose r}\geq{k-o(k)\over \pi_n(\fH){n\choose r}}{n\choose r}=\biggl({1\over\pi(\fH)}-o(1)\biggr)k.
	\]
	
	\textit{Upper bound.} Let $G$ be an $r$-uniform (multi)graph on $n$ vertices with $\ex(G,\fH)<k$, and let $F\subseteq K_n^{(r)}$ be an $\fH$-free subgraph with $e(F)=\ex(K_n^{(r)},\fh)=\pi_n(\fH){n\choose r}$. Let $F'$ be a copy of $F$ chosen uniformly at random from $K_n^{(r)}$ and set $F^*=\{e\in E(G):e\in E(F')\}$, where multiedges are preserved. Certainly as $F$ is $\fH$-free and $\fH$ consists only of simple graphs, $F^*$ is also $\fH$-free. Therefore, since $\pi_n(\fH)\geq\pi(\fH)$,
	\[
	k>\E e(F^*)=\pi_n(\fH)e(G)\geq\pi(\fH)e(G),
	\]
	so $\rk^*_\fh(k)<{k\over\pi(\fH)}$.
\end{proof}

For a family $\fh$ of 2-uniform simple graphs, define $\rho(\fh)=\min_{H\in\fh}\chi(H)$. The Erd\H{o}s--Stone Theorem~\cite{ES46} states that $\pi(\fH)=1-{1\over\rho(\fH)-1}$, so noting that $\pi(\fh)=0$ if and only if $\fh$ contains a bipartite graph, we attain Theorem~\ref{thm:esrank} as a direct corollary.

In the case that $\fh=\{K_t^{(r)}\}$, we conjecture the following:

\begin{conj}\label{conj:cliques}
	If $k=\ex(K_n^{(r)},K_t^{(r)})+1$, then $\rk_{K_t^{(r)}}(k)=\rk_{K_t^{(r)}}^*(k)={n\choose r}$ and the unique extremal graph is $K_n^{(r)}$.
\end{conj}

In Section~\ref{sec:graphs}, we will prove the above conjecture for 2-uniform cliques (see Theorem~\ref{thm:cliques}).

Unfortunately, Theorem~\ref{thm:eshypergraphs} does not give us any useful information about $\rk_\fh(k)$ when $\fh$ contains a degenerate graph. However, we can still say something nontrivial by using the same ideas. We will do so later in Section~\ref{sec:tri}, specifically in relation to 2-uniform graphs (see Theorem~\ref{thm:2degen}).

To end this section, we present a general upper bound on $\rk_H(k)$, which directly follows from the work of Friedgut and Kahn~\cite{FK98} who extended a result of Alon~\cite{alon81}. 

For two hypergraphs $H$ and $G$, let $N(G,H)$ denote the number of copies of $H$ contained in $G$, and let $N(m,H)$ denote the maximum value of $N(G,H)$ taken over all hypergraphs $G$, with $ e(G)=m$. Also, for a hypergraph $H$, we say that $\phi:E(H)\to[0,1]$ is a \emph{fractional cover of $H$} if $\sum_{e\ni v}\phi(e)\geq 1$ for every $v\in V(H)$. The \emph{fractional cover number of $H$}, denoted $\rho^*(H)$ is the minimum value of $\sum_{e\in E(H)}\phi(e)$ where $\phi$ is a fractional cover of $H$.

\begin{theorem}[Friedgut and Kahn~\cite{FK98}]\label{thm:copies}
	For any hypergraph $H$, $N(m,H)=\Theta(m^{\rho^*(H)})$.
\end{theorem}

\begin{prop}\label{prop:genupper}
	Set $\rho^*=\rho^*(H)$ and $s=e(H)$. If $\rho^*<s$, then there is a constant $c=c(H)$ such that 
	\[\rk_H(k)\leq ck^{(s-1)/(s-\rho^*)}.
	\]
\end{prop}
\begin{proof}
	Let $G$ be a graph with $\ex(G,H)<k$ and $ e(G)=m$. By Theorem~\ref{thm:copies}, there is a constant $C=C(H)$ such that $N(G,H)\leq N(m,H)\leq Cm^{\rho^*}$.
	
	We proceed by a standard averaging argument. Let $S\subseteq E(G)$ be a set of edges where each $e\in E(G)$ is included in $S$ independently with probability $p$. Then let $S'\subseteq S$ be attained by removing one edge per copy of $H$ contained in $S$. Thus $S'$ is $H$-free, so
	\[
	k >\E|S'|\geq pm-p^s N(G,H)\geq pm-Cp^sm^{\rho^*}=pm\bigl(1-Cp^{s-1}m^{\rho^*-1}\bigr).
	\]
	Selecting $p^{s-1}m^{\rho^*-1}=1/(sC)$ yields
	\[
	k> \biggl(1-{1\over s}\biggr)\biggl({m^{s-\rho^*}\over sC}\biggr)^{1/(s-1)}.
	\]
	As such, there is some $c=c(H)$ with $
	m<ck^{(s-1)/(s-\rho^*)}$.
\end{proof}

It would be interesting to find a generalization of the fractional cover number so that the above bound would hold for multigraphs as well.

\subsection{Graphs vs.\ Multigraphs}\label{sec:multi}

As mentioned earlier, if $\fH$ is a family of simple graphs, then $\rk_\fH(k)\leq\rk_\fH^*(k)$. In fact, we believe that it would be the case that $\rk_\fh(k)=\rk_\fh^*(k)$ whenever $\fh$ is a family of simple $r$-uniform hypergraphs, as mentioned in Conjecture~\ref{conj:multi}.

This statement appears very difficult to prove in general.
As an aside, $\rk,\rk^*$ have already been actively studied for $\C_o$, the family of odd cycles,
see for example \cite{alon1998bipartite} and 
\cite{bollobas2002better}. Here, results are always phrased as fixing a positive integer $m$ and asking how large a cut (or ``judicious partition'') can be in a (multi)graph with $m$ edges,
or in our terminology, ``How large can $k$ be satisfying $\rk_{\C_o}(k) \leq m$?''
The special relevance to us is that Conjecture 1.2 in \cite{alon1998bipartite} would answer their question for $\rk^*_{\C_o}(k)$, and in particular, show 
$\rk^*_{\C_o}(k)= \rk_{\C_o}(k)$. But as far as we can tell, even this is not known.


While we cannot prove Conjecture~\ref{conj:multi}, we can present the proof of a simple subcase.
\begin{prop}
	Let $\fH$ be a family of simple graphs and $G$ be a multigraph. If each edge of $G$ has the same multiplicity, then there exists a simple graph $G'$ with $e(G')=e(G)$ and $\ex(G',\fH)\leq\ex(G,\fH)$.
\end{prop}
\begin{proof}
	Let $G$ be a multigraph where each edge has multiplicity $r$. Decompose $G$ into simple graphs $G_1,\dots,G_\ell$ and let $G'$ be the disjoint union of $G_1,\dots,G_\ell$, so certainly we have $e(G')=e(G)$. Now, let $F\subseteq G'$ be an $\fH$-free subgraph on $\ex(G',\fH)$ edges and set $F_i=F\cap G_i$. Without loss of generality, suppose $e(F_1)\geq e(F_i)$ for all $i$ and form $F'\subseteq G$ by replacing each edge of $F_1$ by $\ell$ copies. Since $\fH$ consists only of simple graphs, $F'$ is also $\fH$-free, so
	\[
	\ex(G,\fH)\geq e(F')=\ell\cdot e(F_1)\geq \ell\cdot{\ex(G',\fH)\over \ell}=\ex(G',\fH).\qedhere
	\]
\end{proof}

Unfortunately, when $G$ is a multigraph where different edges have different multiplicities, it is unclear whether or not one can construct a simple graph $G'$ with $e(G')= e(G)$ and $\ex(G',\fh)\leq\ex(G,\fh)$. 

Notice that if $\fh$ does not contain a degenerate graph, then $\rk_\fh(k)=\bigl(1+o(1)\bigr)\rk_\fh^*(k)$ by Theorem~\ref{thm:eshypergraphs}. We can also provide the following bound which, unfortunately, is not very strong.

\begin{prop}
	If $\fH$ is a family of simple $r$-uniform graphs, then $\rk_\fh^*(k)\leq\rk_{ k\log k}(\fH)$.
\end{prop}
\begin{proof}
	Both are infinite if $\fh$ contains a sunflower, so we shall suppose that is not the case.
	
	Let $G$ be a multigraph with $\ex(G,\fh)<k$. As above, decompose $G$ into simple graphs $G_1,\dots,G_\ell$ where $G_1\supseteq \dots\supseteq G_\ell$, and let $G'$ be the disjoint union of these graphs, so certainly $e(G')=e(G)$. We now argue that $\ex(G',\fh)<k\log k$, which will establish the claim.
	
	To do this, we first note that as $\ex(G,\fh)<k$, we must have $\ell\leq k-1$ since $\fH$ does not contain $K_r^{(r)}$. Consider any $\fh$-free subgraph $F\subseteq G_i$. As $G_1\supseteq\dots\supseteq G_\ell$, and $\fH$ is a family of simple graphs, we can form an $\fh$-free subgraph $F'\subseteq G$ by replacing every edge of $F$ by $i$ copies. Thus, it must be the case that $\ex(G_i,\fh)<{k\over i}$. As such,
	\[
	\ex(G',\fH)\leq\sum_{i=1}^\ell\ex(G_i,\fH)<\sum_{i=1}^\ell{k\over i}\leq k\log(\ell+1)\leq k\log k.\qedhere
	\]
\end{proof}

Interestingly, an analogue of Conjecture~\ref{conj:multi} fails if we consider non-uniform graphs, and in Section~\ref{sec:nonu}, we give such an example.

\section{2-Uniform Graphs}\label{sec:graphs}

Throughout this section, all graphs will be 2-uniform.

We begin our study of $\rk_\fh(k),\rk_\fh^*(k)$ where $\fh$ is a family of 2-uniform graphs by looking at $\fh=\{K_t\}$. Recall that Theorem~\ref{thm:esrank} told us that $\rk_{K_t}(k),\rk_{K_t}^*(k)=\bigl(1+{1\over t-2}+o(1)\bigr)k$ for all $t\geq 3$. We will find the exact answer for an infinite sequence of values of $k$, thereby proving Conjecture~\ref{conj:cliques} for 2-uniform cliques. To do so, we will rely on a reduction, which will also be used for other results in this section.

We begin by defining a contraction in a graph:

\begin{defn}\label{contract}
If $G$ is a (multi)graph and
$I\subseteq V(G)$,
define $G'=C_{I}(G)$ to be the (multi)graph (possibly with loops) with the same number of edges obtained by contracting together the vertices in $I$. More specifically, write $V(G'):= (V(G) \cup \{z\} )\backslash I$ for some new vertex $z$, and the multiset
$E(G'):= \{C_{I}(e) :e \in E(G)\}$, where
\[
C_{I}(e):= \begin{cases}
zz & \text{ if $e\in{I\choose 2}$;}\\
zx & \text{ if $e=ux$ for some  $u \in I$;}\\
e & \text{ otherwise. }
\end{cases}
\] Here, we think of $C_{I}$ as a bijection between multigraph edge sets. 
\end{defn}

Note that $C_I(G)$ will contain loops whenever $I$ is not an independent set in $G$. In the situations we will use contractions, we will only contract independent sets, but in stating the following definition and lemma, we will allow this to happen.

To apply contractions in determining $\rk_\fh^*(k)$, we provide the following general definition and lemma.

\begin{defn}\label{grelation}
If $\mathcal{G}$ denotes the space of all finite simple graphs and $\mathcal{G}^*$ denotes the space of all finite multigraphs, a function $f:\mathcal{G}^*\to\mathcal{G}$ is called a \emph{graph simplification} if it preserves vertex sets
and containment. That is, for every pair of graphs $G,H$, $V(f(G))=V(G)$ and if $H\subseteq G$, then $f(H)\subseteq f(G)$.
\end{defn}

Examples of graph simplifications include:
\begin{enumerate}
\item $f(G) = G_s$ where $G_s$ is the underlying simple graph of $G$.
\item $ab \in E(f(G)) \Leftrightarrow a,b$ are in the same connected component of $G$,
\item $ab \in E(f(G)) \Leftrightarrow \text{dist}_G(a,b) \leq t$ for some fixed integer $t$.
\end{enumerate}

\begin{lemma}\label{lem:relation}
Let $H$ be a (multi)graph and let $f$ be a graph simplification such that $f(H)$ is a clique for every $H \in \fh$. Suppose that $G$ is a (multi)graph and let $I$ be an independent set in $f(G)$.
If $G'=C_I(G)$, then $\ex(G',\fh)\leq \ex(G,\fh)$.\footnote{We remark that for a general graph simplification $f$, it can be the case that $I$ is an independent set in $f(G)$, but $I$ is \emph{not} an independent set in $G$. As an example, consider the graph simplification defined by $ab\in E(f(G))\Leftrightarrow a,b$ have $t$ vertex-disjoint paths connecting them for $t\geq 2$. In regard to this lemma, if it is the case that $G'$ has loops, then in calculating $\ex(G',\fh)$, we allow using these loops.}
\end{lemma}

\begin{proof}
It suffices to show that if some $F\subseteq G$ contains a copy of $H\in\fh$, then $C_I(F)\subseteq G'$ still contains a copy of some $H'\in\fh$. In fact, more is true; namely, if $H_0\subseteq G$ is a copy of $H$, then $C_I(H_0)\subseteq G'$ contains a copy of $H$. To see this, as $f$ is a graph simplification, $f(H)\simeq f(H_0)\subseteq f(G)$, so as $f(H)$ is a clique, $|I\cap V(H_0)|\leq 1$. In other words, $C_I(H_0)$ is a copy of $H$, possibly with extra multiedges or loops.
\end{proof}

For a graph simplification $f$, we say that $G$ is \emph{$f$-compressed} if $f(G)$ is a clique. Furthermore, we say that $G$ is an \emph{$f$-compressed copy of $G'$} if $G$ is $f$-compressed and there is a sequence of graphs $G'=G_0,G_1,\dots,G_t=G$ such that $G_{i+1}=C_I(G_i)$ for some independent set $I$ in $f(G_i)$. Note that if $G$ is an $f$-compressed copy of $G'$, then $e(G)=e(G')$ (including counting any loops if they exist). With this definition, the following corollary follows immediately from Lemma~\ref{lem:relation}.

\begin{cor}{\label{cor:connect}}
Suppose, as above, that $f$ is a multigraph simplification where $f(H)$ is a clique for every $H\in\fh$. If $G$ is an $f$-compressed copy of $G'$, then $e(G)=e(G')$ and $\ex(G,\fH)\leq\ex(G',\fH)$.
In particular, when computing $\rk_\fh^*(k)$, it suffices to consider graphs $G$ such that $f(G)$ is a clique, i.e.\ 
\[
\rk^*_\fh(k)=\sup \bigl\{e(G): \ex(G,\fh) < k, f(G)\simeq K_{|V(G)|} \bigr\}.
\]
\end{cor}

Before finding the value of $\rk_{K_t}(k)$, we first must recall some properties of Tur\'an graphs. Define $T_{t-1}(n)$ to be the balanced complete $(t-1)$-partite graph on $n$ vertices; Tur\'an's Theorem states that $\ex(K_n,K_t)=e(T_{t-1}(n))$. Additionally, set $\pi_n(t):=\pi_n(K_t)=\ex(K_n,K_t)/{n\choose 2}$, the Tur\'an density of $K_t$. We will use the following observations in the subsequent proof.

\begin{prop}\label{obs:turan}
For all $n\geq 2$, $\pi_n(t)\geq \pi_{n+1}(t)$. Furthermore,
\[
e(T_{t-1}(n))-e(T_{t-1}(n-1))=(n-1)-\bigg\lfloor{n-1\over t-1}\bigg\rfloor.
\]
\end{prop}

We now prove Conjecture~\ref{conj:cliques} for 2-uniform cliques by using an idea of Alon~{\cite[Lemma 2.1]{alon1996bipartite}} in application to chromatic numbers.

\begin{theorem}\label{thm:cliques}
Fix $t\geq 3$ and $n\geq t$. If $k=\ex(K_n,K_t)+1$, then $\rk_{K_t}(k)=\rk_{K_t}^*(k)={n\choose 2}$ and the unique extremal graph for $\rk_{K_t}^*(k)$ (and thus for $\rk_{K_t}(k)$) is $K_n$.
\end{theorem}

\begin{proof} The lower bound is immediate, so we focus only on the upper bound and the classification of the extremal graphs.

Let $G$ be a (multi)graph with $\ex(G,K_t)<k$. Letting $f$ be the ``underlying simple graph'' simplification (example (1) in Definition~\ref{grelation}), as $f(K_t)=K_t$, we may suppose that $G$ is $f$-compressed by Corollary~\ref{cor:connect}. In other words, $G$ is a clique, possibly with parallel edges. Let $n=|V(G)|$ and write $e(G)={r\choose 2}+\ell$ where $0\leq\ell\leq r-1$. Since $G$ is a copy of $K_n$, possibly with parallel edges, we know that $r\geq n$.

Now, let $T$ be a copy of the Tur\'{a}n graph $T_{t-1}(n)$ chosen uniformly at random on $V(G)$, and let $H$ be the multigraph with edge set $\{uv \in E(G): uv\in E(T)\}$ (so that if $u,v$ span multiple edges in $G$ then they either all survive the intersection with $T$ or all do not). Since any such $H$ is $K_t$-free, we calculate
\begin{align*}
	\ex(G,K_t) &\geq \E e(H) = e(G)\cdot \pi_n(t)\geq e(G)\cdot \pi_r(t)\\
	&=\biggl({r\choose 2}+\ell\biggr)\cdot{\ex(K_r,K_t)\over {r\choose 2}}\\
	&=\ex(K_r,K_t)+\ell\cdot\pi_r(t).
\end{align*}

However, by assumption, $\ex(G,K_t)<k=\ex(K_n,K_t)+1$, so we have $\ex(K_r,K_t)+\ell\cdot\pi_r(t)\leq \ex(K_n,K_t)$. However, since $r\geq n$, $\pi_r(t)>0$ and $\ell\geq 0$, the only way for this to happen is if $\ell=0$ and $n=r$; hence, $G\simeq K_n$. In particular, $\rk_{K_t}(k)\leq\rk_{K_t}^*(k)\leq{n\choose 2}$.

Now, suppose that $G$ is any (multi)graph on ${n\choose 2}$ edges with $\ex(G,K_t)<k$. Let $G=G_0,G_1,\dots,G_q=G^*$ where $G^*$ is $f$-compressed and $G_{i+1}=C_{xy}(G_i)$ for some $xy\notin E(G_i)$. By the statement above, we know that $G^*\simeq K_n$. 

Suppose that $G\not\simeq K_n$, so in particular $q\geq 1$. Let $u,v\in V(G_{q-1})$ be such that $G^*=C_{uv}(G_{q-1})$. For ease of notation, we will write $N(x)=N_{G_{q-1}}(x)$ for the remainder of the proof.

Since $uv\notin E(G_{q-1})$, we know that $G_{q-1}\not\simeq K_n$. However, $K_n\simeq G^*=C_{uv}(G_{q-1})$, so $|V(G_{q-1}|=n+1$. Also, we must have $N(u)\cup N(v)=V(G_{q-1})\setminus\{u,v\}$ and $V(G_{q-1})\setminus\{u,v\}$ must induce a copy of $K_{r-1}$. Furthermore, $G^*$ is a simple graph, so it must be the case that $G_{q-1}$ is simple as well; moreover $N(u)\cap N(v)=\varnothing$ otherwise $G^*$ would contain a multiedge upon contracting $uv$. We check that such a graph has a $K_t$-free subgraph which is too large.

Indeed, first suppose $|N(u)|<\bigl\lfloor{n-1\over t-1}\bigr\rfloor$. Then let $T$ be a copy of $T_{t-1}(n-1)$ contained in $V(G_{q-1})\setminus\{u,v\}$ with parts $X_1,\dots,X_{t-1}$ where $X_1\supseteq N(u)$ and $|X_1|=\bigl\lfloor{n-1\over t-1}\bigr\rfloor$. Then if $H$ is the subgraph consisting of the edges in $T$ along with the edges incident to $u$ and edges of the form $\{vx:x\notin X_{t-1}\}$, we find that $H\subseteq T_{t-1}(n+1)$ as $uv\notin E(G_{q-1})$, so $H$ is $K_t$-free. Additionally,
\begin{align*}
e(H) &= e(T_{t-1}(n-1))+|N(u)|+|N(v)|-|X_1\cap N(v)|\\
&\geq e(T_{t-1}(n-1))+(n-1)-\biggl(\biggl\lfloor{n-1\over t-1}\biggr\rfloor-1\biggr)\\
&= e(T_{t-1}(n))+1=k,
\end{align*}
a contradiction. Thus, we may suppose that $|N(u)|,|N(v)|\geq\bigl\lfloor{n-1\over t-1}\bigr\rfloor$. Additionally, since $|N(u)|+|N(v)|=n-1$, we have, without loss of generality, $|N(v)|\geq\bigl\lceil{n-1\over t-1}\bigr\rceil$. As such, let $T$ be a copy of $T_{t-1}(n-1)$ contained in $V(G_{q-1})\setminus\{u,v\}$ with parts $X_1,\dots,X_{t-1}$ where $X_1\subseteq N(u)$ and $X_2\subseteq N(v)$. Now, let $H$ consist of $T$ along with all edges incident to $u$ or $v$. Since $uv\notin E(G_{q-1})$, $H$ is again a subgraph of $T_{t-1}(n+1)$, and so is $K_t$-free. However,
\begin{align*}
e(H) &= e(T_{t-1}(n-1))+|N(u)|+|N(v)|\\
&= e(T_{t-1}(n-1))+n-1\\
&= e(T_{t-1}(n))+\biggl\lfloor{n-1\over t-1}\biggr\rfloor\geq k,
\end{align*}
since $n\geq t$; another contradiction. We conclude that any (multi)graph $G$ with $e(G)={n\choose 2}$ and $\ex(G,K_t)<k$ must be $K_n$.
\end{proof}

It is not clear what the precise value of $\rk_{K_t}(k)$ and $\rk_{K_t}^*(k)$ are when $k\neq\ex(K_n,K_t)+1$ for any $n$, but we conjecture the following:

\begin{conj}
For positive integers $r_1\geq\dots\geq r_\ell$, let $K(r_1,\dots,r_\ell)$ be the multigraph consisting of ``nested'' copies of $K_{r_i}$: that is, on vertex set $[r_1]$, we overlay a copy of $K_{r_i}$ on $[r_i]$ for every $i$ (thus, the maximum edge-weight is $\ell$, provided every $r_i \geq 2$). Then for every $k$, there exist positive integers $r_1\geq\dots\geq r_\ell$ such that $K(r_1,\dots,r_\ell)$ is extremal for $\rk_{K_t}^*(k)$.
\end{conj}

We do not even have a conjectured extremal example for $\rk_{K_t}(k)$.

\subsection{Cycles}

We begin this section with a simple result related to the graphic matroid rank of a graph.

\begin{theorem}\label{thm:cycles}
If $\C:=\{C_3,C_4,\dots\}$ is the set of all cycles, then $\rk_\C(k)=\rk_\C^*(k) = \binom{k}{2}$. Furthermore, the only extremal graph for $\rk_\C^*(k)$ (and thus for $\rk_\C(k)$) is $K_k$.
\end{theorem}

\begin{proof}
Note that any $k$-edge subgraph of $K_k$ contains a cycle; hence $\rk_\C(k) \geq \binom{k}{2}$.

Now, suppose that $G$ is some \emph{connected} (multi)graph with $\ex(G,\C)<k$. If $G$ has at least $k+1$ vertices, then $G$ has a spanning tree with at least $k$ edges; a contradiction. Hence, $|V(G)|\leq k$.

Additionally, for any fixed $v\in V(G)$, the set of edges incident to $v$ is $\C$-free, so $\Delta(G)\leq k-1$. Thus,
\[
e(G)={1\over 2}\sum_{v\in V(G)}d(v)\leq{1\over 2}|V(G)|\Delta(G)\leq{k\choose 2}.
\]
If we have equality, then certainly $|V(G)|=k$. However, $G$ cannot have any edge of multiplicity two or higher, otherwise, extending this multiedge to a spanning tree of $G$ (including the multiedges), we will find a $\C$-free subgraph of $G$ with at least $k$ edges. Thus, $G$ must be simple, so we have $e(G)\leq {k\choose 2}$ with equality if and only if $G\simeq K_k$. Since every cycle is connected, we have proved the upper bound by taking $f$ to be the connectedness graph simplification (example (2) in Definition~\ref{grelation}) in Corollary~\ref{cor:connect}.

We now wish to argue that if $e(G)={k\choose 2}$ and $\ex(G,\C)<k$, then $G\simeq K_k$. If $G$ is connected, we already know this to be the case, so suppose $G$ is disconnected. Letting $I\subseteq V(G)$ consist of one vertex from each connected component of $G$, we see that $C_I(G)\simeq K_k$ by the connected case and Corollary~\ref{cor:connect}. However, this implies that $K_k$ has a cut-vertex, which is not true. Thus $G$ must have been connected in the first place so $G\simeq K_k$.
\end{proof}

If $\fh$ does not contain a bipartite graph, then the asymptotic value of $\rk_{\fh}(k)$ is determined by Theorem~\ref{thm:esrank}. In contrast, it is natural to ask about this extremal function for the family of all even cycles $\C_e=\{C_4,C_6,\ldots\}$. We now prove Theorem~\ref{thm:evencycles}, which stated:

For $k\geq 4$, $\rk_{\C_e}(k)=\bigl\lfloor{k^2\over 4}\bigr\rfloor$. Furthermore, the only extremal graph for $\rk_{\C_e}(k)$ is the balanced complete bipartite graph on $k$ vertices, unless $k =5$.


Along the way, we will need a classical strengthening of the Mantel bound:

\begin{lemma}\label{lemma:mantelind}
Any triangle-free graph on $n$ vertices with independence number $\alpha$ has $\leq \alpha (n-\alpha)$ edges.
\end{lemma}

\begin{proof}[Proof of Theorem~\ref{thm:evencycles}]
\textit{Lower bound.} Let $G$ be the balanced complete bipartite graph on $k$ vertices. Naturally, any $k$ edges from $G$ contain a cycle, necessarily even as $G$ is bipartite. Hence, $\rk_{\C_e}(k)\geq e(G)=\bigl\lfloor{k\over 2}\bigr\rfloor\bigl\lceil{k\over 2}\bigr\rceil=\bigl\lfloor{k^2\over 4}\bigr\rfloor$.

For the upper bound, we again look to use Corollary~\ref{cor:connect}, and first prove the connected case.

\begin{lemma}
If $G$ is \emph{connected} with $ e(G) 
\geq \bigl\lfloor \frac{k^2}{4} \bigr\rfloor$, then $\ex(G,\C_e) \geq k-1$, with equality if and only if $G \simeq K_{\lceil k/2\rceil, \lfloor k/2\rfloor}$ or, in the case of $k=5$, $G\simeq K_4$.  
\end{lemma}

\begin{proof}
Let $G$ be any connected graph on $n$ vertices with $\ex(G,\C_e)\leq k-1$.
For any spanning tree $F$ of $G$, $F$ contains no even cycle, so $e(F)\leq k-1$, or in other words, $n\leq k$. As such,
set $k=n+q$, and assume 
\[
e(G) \geq \biggl\lfloor {k^2\over 4}\biggr\rfloor={n^2\over 4}+{2nq+q^2-\mathbf{1}_{\text{k odd}}\over 4} \geq {n^2\over 4}+\frac{2nq + q^2 -1}{4},
\]
but that $G$ is not the complete balanced bipartite graph. Then as $n \leq k$, we know, by the uniqueness of the Tur\'{a}n graph, that $G$ contains a triangle. We will attempt to use the triangles in $G$ to build a large $\C_e$-free subgraph.

Say $T \subseteq G$ is a ``triangle forest'' with $t$ triangles if $E(T)$ is a collection of $t$ edge-disjoint triangles such that the removal of any one edge from each triangle forms a forest. More formally, $T$ is a triangle forest if every vertex of $T$ is contained in some triangle and each 2-connected block of $T$ is a triangle (of which there are $t$). In particular, the only cycles within such a $T$ are the $t$ triangles. So we may extend $T$ to a spanning subgraph $H$ (using connectivity) with no additional cycles, thus $H$ is still $\C_e$-free. We deduce that
$(n-1)+t=e(H) \leq k-1$. Thus $t \leq q$ for any such $T$. Since $G$ has at least one triangle, we note $q\geq 1$.

 Now, take such a triangle forest $T$ with:
\begin{enumerate}
\item $e(T)$ (and hence $t$) as large as possible,
\item Subject to (1), if $T=T_1 \cup \dots\cup T_\ell$ is a decomposition of $T$ into connected components where $|T_1| \geq \dots \geq |T_\ell|$, then $(|T_1|, \dots, |T_\ell|) $ is maximal in the lexicographic ordering.
\end{enumerate}
By the lexicographic order, we mean that $(a_1, \dots, a_\ell) \succ (b_1, \dots, b_{\ell'}) \Leftrightarrow a_j > b_j$ for $j:= \min\{i: a_i \neq b_i\}$.
Such a lexicographic maximal $T$ means there is no $v \in T_i$ with 2 edges to the same triangle in $T_j$ for any $i < j$. If this were not the case and $wxy$ was the triangle with both $vx,vy \in E(G)$, then let $T':= (T \cup \{vx,vy\}) \backslash \{ wx,wy \}$ (see Figure~\ref{fig:triangle}). $T'$ is a triangle forest with the same number of edges as $T$, with $|T_j'|=|T_j|$ for all $j<i$ yet $|T_i'|\geq|T_i|+2$, so $T'$ is lexicographically larger than $T$, contradicting (2).

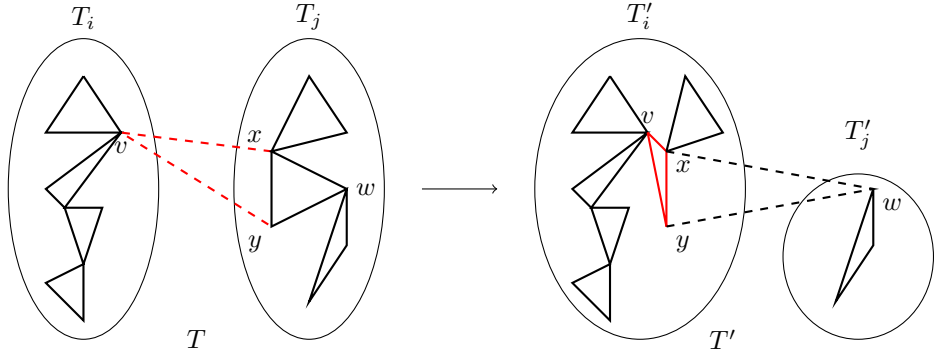
\begin{figure}[h]
\centering
\begin{tikzpicture}

\draw (0,0) circle [x radius = 1, y radius = 2];
\draw (3,0) circle [x radius = 1, y radius = 2];
\node at (1.5,-2) {$T$};
\node [above] at (0,2) {$T_i$};
\node [above] at (3,2) {$T_j$};

\draw [thick] (0,1.5)--(-0.5,0.75)--(0.5,0.75)
--(-0.5,0)--(-0.25,-0.25)
--(0,-1)--(-0.5,-1.25)--(0,-1.75)--(0,-1)--(0.25,-0.25)--(-0.25,-0.25)--(0.5,0.75)--(0,1.5);

\draw [thick, red, dashed] (2.5,0.5)--(0.5,0.75)--(2.5,-0.5);

\draw [thick] (2.5,0.5)--(2.5,-0.5)--(3.5,0)--(3,-1.5)--(3.5,-0.75)--(3.5,0)--(2.5,0.5)--(3.5,0.75)--(3,1.5)--(2.5,0.5);

\node [right] at (3.5,0) {$w$};
\node [above left] at (2.5,0.5) {$x$}; 
\node [below left] at (2.5,-0.5) {$y$}; 
\node [below] at (0.5,0.75) {$v$}; 

\draw [->] (4.5,0)--(5.5,0);

\draw (7.4,0) circle [x radius = 1.4, y radius = 2];
\draw (10.3,-0.9) circle [x radius = 1, y radius = 1.1];
\node at (8.5,-2) {$T'$};
\node [above] at (7.4,2) {$T_i'$};
\node [above] at (10.3,.4) {$T_j'$};

\draw [thick] (7,1.5)--(6.5,0.75)--(7.5,0.75)
--(6.5,0)--(6.75,-0.25)
--(7,-1)--(6.5,-1.25)--(7,-1.75)--(7,-1)--(7.25,-0.25)--(6.75,-0.25)--(7.5,0.75)--(7,1.5);
\node [below right] at (10.5,0) {$w$};
\node [below right] at (7.75,0.5) {$x$}; 
\node [below right] at (7.75,-0.5) {$y$}; 
\node [above] at (7.5,0.75) {$v$}; 

\draw [red, thick] (7.5,0.75)--(7.75,0.5)--(7.75,-0.5)--(7.5,0.75);
\draw [thick] (7.75,0.5)--(8.5,0.75)--(8,1.5)--(7.75,0.5);
\draw [thick] (10.5,0)--(10,-1.5)--(10.5,-0.75)--(10.5,0);

\draw [dashed, thick] (7.75,0.5)--(10.5,0)--(7.75,-0.5); 

\end{tikzpicture}
\caption{Finding a lexicographically larger triangle forest in the case where some vertex in $T_i$ has two edges to the same triangle in $T_j$. \label{fig:triangle}}
\end{figure}

Thus, if $T_j$ consists of $t_j$ triangles for every $j$ (so that $|T_j|=2t_j +1$ and $t= \sum_j t_j$), then whenever $i < j$, every $v \in T_i$  has at most $t_j$ edges to $T_j$.
 Summing over all $v \in T_i$ gives $e[T_i,T_j] \leq (2t_i+1)t_j$.

We now attempt to bound the remaining edges in $G$.
Crudely, $e[T_i] \leq \binom{|T_i|}{2} =2t_i^2 + t_i$ for all $i$.

\emph{Case 1:} $|V(T)| \leq \frac{n}{2}$.\\
Let $G':= G \setminus \bigcup_i G[T_i]$. As $T$ is maximal, $G'$ must be triangle-free,
so certainly $e(G') \leq \frac{n^2}{4}.$ Therefore,

\[
\frac{n^2}{4}+ \frac{2nq+q^2-1}{4}
\leq  e(G)= e(T)+e(G')
\leq \sum_{i=1}^{\ell} (2t_i^2 + t_i) + \frac{n^2}{4},
\]
and so
\[
\frac{2nq+q^2-1}{4} \leq t_1 \sum_{i=1}^{\ell} (2t_i +1)
=t_1 |V(T)| \leq t |V(T)| \leq q \cdot \frac{n}{2}.
\]
Thus, $q=1$ as we supposed that $q\geq 1$, so we have equality everywhere. In particular, $t_1=t=q=1$, so $T$ is a single triangle, $|V(T)|=\frac{n}{2} \Rightarrow n=6$, and $e(G')=\frac{n^2}{4} \Rightarrow G'=K_{3,3}$. Since $G$ is therefore a 6-vertex, edge-disjoint union of $K_{3,3}$ with a triangle, this uniquely determines $G$ as $K_6 \backslash K_3$, and this $G$ still has $\ex(K_6 \backslash K_3,\C_e)\geq 7=n+q=k$ (see Figure~\ref{fig:k6k3}); a contradiction.

\begin{figure}[h]
\centering
\begin{tikzpicture}
\draw [red, ultra thick] (0,0)--(0,2)--(-1,2)--(-1,0)--(0,2)--(1,0)--(1,2)--(0,2);
\draw [thick] (1,0)--(-1,2)--(0,0)--(1,2)--(-1,0);
\draw [thick] (-1,2) to [out=30,in=150] (1,2);
\draw [fill] (-1,2) circle  [radius=0.08];
\draw [fill] (0,2) circle [radius=0.08];
\draw [fill] (1,2) circle [radius=0.08];
\draw [fill] (-1,0) circle [radius=0.08];
\draw [fill] (0,0) circle [radius=0.08];
\draw [fill] (1,0) circle [radius=0.08];
\end{tikzpicture}
\caption{A $\C_e$-free subgraph of $K_6\backslash K_3$ with 7 edges.\label{fig:k6k3}}
\end{figure}

\emph{Case 2:} $|V(T)|> \frac{n}{2}$.\\
In this case, for the triangle-free graph $G'':=G \backslash G[T]
=G' \backslash \bigcup_{i,j} G[T_i,T_j]$, $V(T)$ spans an independent set in $G''$, so
$\alpha(G'') \geq|V(T)|=2t+\ell>{n\over 2}$. Applying Lemma \ref{lemma:mantelind}, and
noting $x(n-x)$ is strictly decreasing for $x\geq n/2$, we have $e(G'')\leq \alpha(G'')(n-\alpha(G'')) \leq (2t+\ell)(n-(2t+\ell))$.

We run a similar calculation in this case:

\begin{align*}
\bigg\lfloor \frac{k^2}{4}\bigg\rfloor \leq  e(G) & \leq \sum_{i = 1}^{\ell} e[T_i]
+ \sum_{i < j} e[T_i,T_j] + e(G'')\\
& \leq \sum_{i=1}^{\ell} \bigl(2t_i^2 +t_i\bigr) + \sum_{i<j}\bigl((2t_i+1)t_j\bigr)
+\bigl(n-(2t+\ell)\bigr)(2t+\ell)\\
& \leq \sum_{i=1}^{\ell} (2t_i^2+t_i) +
\sum_{i \neq j}\biggl( t_it_j + \frac{t_i + t_j}{4}\biggr)
+n(2t+\ell)-(2t+\ell)^2\\
& =t^2 + \sum_{i=1}^{\ell} t_i^2 + 
\biggl( \frac{\ell +1}{2} \biggr) t + n(2t+\ell)-(2t+\ell)^2,
\end{align*}
so
\begin{align*}
n^2+2qn+q^2 - 4n(2t+\ell)+4(2t+\ell)^2 - \mathbf{1}_{k \text{ odd}}
& \leq 8t^2 + 2 (\ell+1) t\\
\Rightarrow
\bigl(n+q-2(2t+\ell)\bigr)^2 +4q(2t+\ell)-\mathbf{1}_{k \text{ odd}}
& \leq 8t^2 +2(\ell+1)t \leq 8qt + 2(2\ell )q.\\
\end{align*}
It follows $|k-2(2t + \ell)| \leq \mathbf{1}_{k \text{ odd}}$. But the reverse is true whether $k$ is even or odd, hence we again obtain all inequalities above at equality. So certainly $t=q$, $\ell=1$, $G[V(T)]=G[V(T_1)]$ is a clique, and
$\alpha(G'')=2t+\ell$, so
$e(G'')= \sum_{v \notin V(T)} \deg_{G''}(v)=(n-(2t+\ell))(2t+\ell)$. As such, $G''[\overline{V(T)}]$ is empty and 
$\deg_{G''}(v)=2t+\ell$ for every $v \notin V(T)$, so $G''$ is the complete bipartite graph on $[V(T),\overline{V(T)}]$. Putting this together with the clique on $V(T)$, deduce $G \simeq K_n \backslash K_r$, where $r=n-(2t+\ell)=n-(2q+1)$.

We know $k-(4q+2)=:\epsilon \in \{0, \pm 1 \}$, so
$r=(k-q)-(2q+1)=q+1+\epsilon$. Now, if $r\geq q+1$, we can find a triangle forest $F$ with $q+1$ triangles (contradicting maximality of $T$ as $t\leq q$) by taking a path on $2(q+1)$ edges with $q+1$ vertices in $\overline{V(T)}$ and $q+2 \leq 2q+1$ among $V(T)$, and completing the $q+1$ edge-disjoint copies of $P_3$ into $K_3$'s using the edges from inside $T$ (See Figure~\ref{fig:triforest}).
Furthermore, if $q \geq 2,$ then we may similarly choose $P$ by instead taking $q \leq q+1+\epsilon$ vertices of $\overline{V(T)}$ and $ q+3 \leq 2q+1$ vertices of $V(T)$.
Otherwise, $\epsilon = -1$ and $q=1$. In this case, we deduce that $G\simeq K_4$, which does have $\ex(K_4,\C_e)=\ex(K_{2,3},\C_e)=4$.
\end{proof}

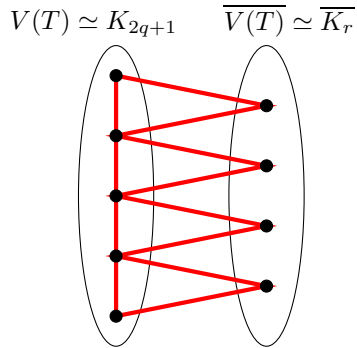
\begin{figure}[ht]
\centering
\begin{tikzpicture}
\draw (0,0) circle [x radius = .5, y radius = 2];
\draw (2,0) circle [x radius =.5, y radius=2];
\node [above] at (-.3,2) {$V(T)\simeq K_{2q+1}$};
\node [above] at (2.3,2) {$\overline{V(T)}\simeq\overline{K_r}$};
\draw[red,ultra thick] (0,1.6)--(2,1.2)--(0,.8)--(2,.4)--(0,0)--(2,-.4)--(0,-.8)--(2,-1.2)--(0,-1.6)
(0,1.6)--(0,-1.6);
\draw [fill] (0,1.6) circle  [radius=0.08];
\draw [fill] (2,1.2) circle  [radius=0.08];
\draw [fill] (0,.8) circle  [radius=0.08];
\draw [fill] (2,.4) circle  [radius=0.08];
\draw [fill] (0,0) circle  [radius=0.08];
\draw [fill] (2,-.4) circle  [radius=0.08];
\draw [fill] (0,-.8) circle  [radius=0.08];
\draw [fill] (2,-1.2) circle  [radius=0.08];
\draw [fill] (0,-1.6) circle  [radius=0.08];
\end{tikzpicture}
\caption{\label{fig:triforest} A large triangle forest contained in $K_n\setminus K_r$.}
\end{figure}

\textit{Upper bound.}
If $G$ is now arbitrary with $ e(G)\geq \bigl\lfloor \frac{k^2}{4} \bigr\rfloor$, then forming any $I\subseteq V(G)$
 with one vertex from each connected component gives $\ex(G,\C_e) \geq \ex(C_{I}(G),\C_e) \geq k-1$.
If we have equality here, we know $C_I(G)$ is necessarily $K_{\lfloor k/2 \rfloor, \lceil k/2\rceil }$ (or $K_4$) by the lemma, 
yet none of these graphs have a cut-vertex for $k \geq 4$. Hence, $G$ must have been connected in the first place, so $G$ is one of the claimed extremal graphs.
\end{proof}

Sadly, the above argument is very specific to simple graphs, so we have been unable to determine $\rk^*_{\C_e}(k)$.
\begin{conj}
	For $k\geq 6$, $\rk_{\C_e}^*(k)=\left\lfloor \frac{k^2}{4} \right\rfloor$, and $K_{\lfloor k/2\rfloor,\lceil k/2\rceil}$ is also the unique extremal graph for $\rk_{\C_e}^*(k)$.
\end{conj}

\subsection{Small Graphs}

In this section, we will explore $\rk_\fh(k)$ where $\fh$ is a collection of small graphs. At the end of this section, we also give a complete classification of the families which have $\rk_\fh(k)=\infty$. Throughout this section, we will only focus on simple host graphs.

Recall that $P_t$ denotes the path on $t$ \emph{edges}.

We first turn our attention to determining $\rk_{P_3}(k)$. We note that $H$ is $P_3$-free if and only if $H$ is the vertex-disjoint union of triangles, stars and isolated vertices. The following graphs arise naturally in determining $\rk_{P_3}(k)$ and classify most extremal graphs.

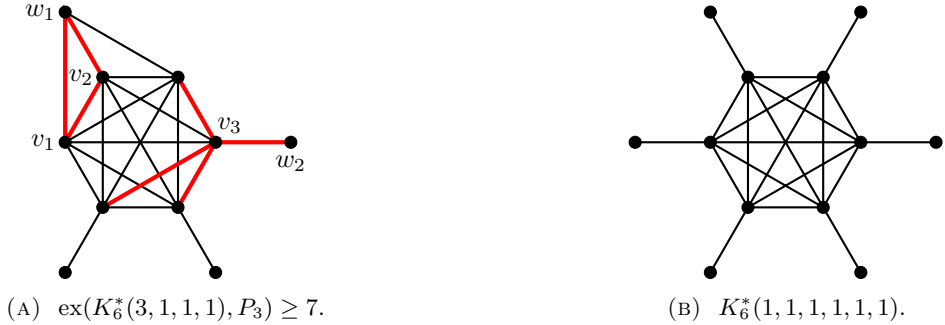
\begin{figure}[ht]
\centering
\begin{subfigure}[t]{0.5\textwidth}
        \centering
        \begin{tikzpicture}
\draw[thick] (-0.5,-0.866)--(0,0)--(1,0)--(1.5,0.866)--(1,1.732)
--(0,1.732)--(-0.5,0.866)--(0,0)
--(1.5,0.866)--(0,1.732)--(0,0)--(1,1.732)--(-0.5,0.866)--(1,0)--(1,1.732)--(-0.5,2.598);
\draw[thick] (-0.5,0.866)--(1.5,0.866)--(2.5,0.866);
\draw[thick] (1.5,-0.866)--(1,0)--(0,1.732)--(-0.5,2.598);
\draw [red,ultra thick] (-0.5,0.866)--(0,1.732)--(-0.5,2.598)--(-0.5,0.866);
\draw[red,ultra thick] (1,1.732)--(1.5,0.866)--(2.5,0.866) ;
\draw[red,ultra thick] (1,0)--(1.5,0.866)--(0,0) ;
\node [left] at (-0.5,2.598) {
$w_1$
};
\node [below] at (2.5,0.806) {
$w_2$
};
\node [above right] at (1.4,0.866) {
$v_3$
};
\node [left] at (-0.5, 0.866) {
$v_1$
};
\node [left] at (0, 1.732) {
$v_2$
};

\draw[fill] (0,0) circle [radius=0.08];
\draw[fill] (1,0) circle [radius=0.08];
\draw[fill] (1.5,0.866) circle [radius=0.08];
\draw[fill] (1,1.732) circle [radius=0.08];
\draw[fill] (0,1.732) circle [radius=0.08];
\draw[fill] (-0.5,0.866) circle [radius=0.08];

\draw[fill] (-0.5,2.598) circle [radius=0.08];
\draw[fill] (-0.5,-0.866) circle [radius=0.08];
\draw[fill] (1.5,-0.866) circle [radius=0.08];
\draw[fill] (2.5,0.866) circle [radius=0.08];
        \end{tikzpicture}
        \caption{\label{fig:pendenta} $\ex(K_6^*(3,1,1,1),P_3)\geq 7$.}
    \end{subfigure}%
\begin{subfigure}[t]{0.5\textwidth}
\begin{tikzpicture}
\node at (1.5,0) {};
\draw[thick] (4.5,-0.866)--(5,0)--(6,0)--(6.5,0.866)--(6,1.732)
--(5,1.732)--(4.5,0.866)--(5,0)
--(6.5,0.866)--(5,1.732)--(5,0)--(6,1.732)--(4.5,0.866)--(6,0)--(6,1.732)--(6.5,2.598);
\draw[thick] (3.5,0.866)--(4.5,0.866)--(6.5,0.866)--(7.5,0.866);
\draw[thick] (6.5,-0.866)--(6,0)--(5,1.732)--(4.5,2.598);
\node [right] at (6,1.732) {
};
\node [left] at (6.5,2.598) {
};

\draw[fill] (5,0) circle [radius=0.08];
\draw[fill] (6,0) circle [radius=0.08];
\draw[fill] (6.5,0.866) circle [radius=0.08];
\draw[fill] (6,1.732) circle [radius=0.08];
\draw[fill] (5,1.732) circle [radius=0.08];
\draw[fill] (4.5,0.866) circle [radius=0.08];

\draw[fill] (4.5,2.598) circle [radius=0.08];
\draw[fill] (3.5,0.866) circle [radius=0.08];
\draw[fill] (4.5,-0.866) circle [radius=0.08];
\draw[fill] (6.5,-0.866) circle [radius=0.08];
\draw[fill] (7.5,0.866) circle [radius=0.08];
\draw[fill] (6.5,2.598) circle [radius=0.08];
\end{tikzpicture}
\caption{\label{fig:pendentb} $K_6^*(1,1,1,1,1,1)$.}
\end{subfigure}
\caption{Pendant graphs. \label{fig:pendent}}
\end{figure}

\begin{defn}
For fixed positive integers $k, r_1, r_2, \dots, r_s$ with $\sum_{i=1}^s r_i=k$, 
define the \emph{pendant} graph $K_k^*(r_1,\dots, r_s)$ as follows.
Take a clique on some $k$-vertex set $\{v_1, \dots, v_k\}$, called the \emph{core}, and additional vertices $\{w_1,\dots, w_s\}$, called the pendants. Partition $\{v_1,\dots,v_k\}=W_1\cup\dots\cup W_s$ where $|W_i|=r_i$ and connect $w_i$ to the vertices in $W_i$. See Figure~\ref{fig:pendent}.

As such, the degree sequence of $K_k^*(r_1,\dots,r_s)$ is $(\underbrace{k,\dots,k}_{k},r_1,\dots,r_s)$ and $e(K_k^*(r_1,\dots,r_s))={k+1\choose 2}$. 
\end{defn}

\begin{lemma}\label{lem:p3cons}
Let $k\geq 4$ and let $r_1,\dots,r_s$ be positive integers with $\sum_{i=1}^s r_i=k-1$. We have\\ $\ex(K_{k-1}^*(r_1,\dots,r_s),P_3)\geq k-1$, where equality holds if either $r_i=1$ for all $i$, or $3\nmid k$ and $r_1=k-1$. Thus, $\rk_{P_3}(k)\geq{k\choose 2}$. (Though unnecessary, this ``if'' is also an ``only if''.)
\end{lemma}

\begin{proof}
Every vertex in the core of $G:=K_{k-1}^*(r_1,\dots,r_s)$ has degree $k-1$, so $\ex(G,P_3)\geq k-1$ is immediate by taking any star centered at a core vertex of $G$.

Now, if $r_1=k-1$, then $G\simeq K_k$, and it is well-known that $\ex(K_k,P_3)=k-1$ if $3\nmid k$. If $(r_1,\dots,r_s)=(1,\dots,1)$, then let $U$ denote the core of $G$.
Now let $H \subseteq G$ be any $P_3$-free
subgraph, so $H$ is a vertex-disjoint union of triangles, stars and isolated vertices. We wish to show $e(H) \leq k-1$.
Now, no triangle $T$ in $H$ can contain a pendant vertex, so each $V(T) \subseteq U$,
and every star contains at most one; hence
$|V(S) \cap U| \geq |V(S)|-1$ for each star $S$. Hence, splitting up $H$ into components:
\begin{align*}
e(H) &= \sum_{\substack{T \subseteq H \\ T \text{ triangle} } } e(T) +
\sum_{\substack{S \subseteq H \\ S \text{ star} } } e(S)\\
&= \sum_{\substack{T \subseteq H \\ T \text{ triangle} } } |V(T)| +
\sum_{\substack{S \subset H \\ S \text{ star} } } \bigl(|V(S)|-1\bigr) \\
& \leq 
\sum_{\substack{T \subseteq H \\ T \text{ triangle} } } |V(T) \cap U| +
\sum_{\substack{S \subseteq H \\ S \text{ star} } } |V(S) \cap U| \leq |U|=k-1. \qedhere 
\end{align*}

\end{proof}

Together with $\ex(C_4,P_3)=3$, this establishes $\rk_{P_3}(k) \geq \binom{k}{2} + \mathbf{1}_{k=3}$. We now show this is the correct value.

\begin{theorem}\label{thm:p3}
If $G$ is a graph with $\ex(G,P_3)< k$,
then $ e(G) \leq \binom{k}{2}+\mathbf{1}_{k=3}$. Thus, $\rk_
{P_3}(k)=\binom{k}{2}+\mathbf{1}_{k=3}$.
\end{theorem}

\begin{proof}
The instances of $k \leq 3$ can be checked by hand, so we assume $k \geq 4$ and proceed by strong induction on $k$.
Note that $\ex(G,P_3) \geq \Delta:=\Delta(G)$, so
$\Delta \leq k-1$.

Firstly, suppose $G$ contains a triangle $T=xyz$. If $H\subseteq G[V\setminus T]=:G'$ is $P_3$-free, then $H\cup T$ is also $P_3$-free, so $\ex(G',P_3)<k-3$. Thus, by induction, $e(G')\leq{k-3\choose 2}+\mathbf{1}_{k-3=3}$. Now, as $\Delta\leq k-1$, $x,y,z$ all have at most  $k-3$ neighbors outside $T$, so 
\[
 e(G)\leq e[V\setminus T]+3(k-3)+3\leq{k-3\choose 2}+\mathbf{1}_{k-3=3}+3k-6={k\choose 2}+\mathbf{1}_{k-3=3},
\]

This is the desired bound unless $k=6$.
Here, one may establish by hand that any graph $G$ with $\binom{k}{2}+1=16$ edges, a copy of $T \cup C_4$, and 3 additional edges off each vertex of $T$, still has $\ex(G,P_3) \geq 6$.

We now suppose that $G$ is triangle-free. Fix $xy\in E(G)$, then this means $N(x)\cap N(y)=\varnothing$. Taking maximal stars with centers $x$ and $y$ (except for the edge $xy$), yields a $P_3$-free subgraph of $G$, so $k>\ex(G,P_3)\geq(\deg(x)-1)+(\deg(y)-1)$. We have thus shown $\deg(x)+\deg(y)\leq k+1$ for every edge $xy$.

Now say there is some edge $xy$ with $\deg(x)+\deg(y)\leq k$. Set $G':=G\setminus\{x,y\}$. Here, adding edge $xy$ to any $P_3$-free subgraph of $G'$ shows that $\ex(G',P_3)\leq\ex(G,P_3)-1<k-1$. Thus, induction gives $e(G) \leq e(G') + (k-1) \leq \binom{k-1}{2} + (k-1) = \binom{k}{2}+\mathbf{1}_{k-1=3}$, and the $k=4$ case may be manually strengthened to remove the extra $+\mathbf{1}_{k-1=3}$.

Hence, we may suppose $\deg(x)+\deg(y)=k+1$ for every $xy\in E(G)$. Fix $x$ and suppose first that $d:=\deg(x)\neq{k+1\over 2}$. Letting $C$ denote the connected component of $G$ containing $x$, we can partition $C=A\cup B$ where $A=\{u:\deg(u)=d\}$ and $B=\{u:\deg(u)=k+1-d\}$. As $\deg(u)+\deg(v)=k+1$ for every $uv\in E(G)$ and $d\neq{k+1\over 2}$, $G[C]$ is a bipartite graph with parts $A,B$. Now, for any $u\in A$ and $v\in B$, by considering stars centered at $u$ and $v$ (except for the edge $uv$ if it exists), we find
\[
k>\ex(G,P_3)\geq\ex(G[C],P_3)+\ex(G[V\setminus C],P_3)\geq|N(u)\setminus \{v\}|+|N(v)\setminus\{u\}|=k+1-2\cdot\mathbf{1}_{uv\in E(G)}.
\]
From this, we immediately find that $G[V\setminus C]$ is empty, and as the above holds for any $u,v$, we know that $G[C]$ is a complete bipartite graph. Furthermore, since $C$ is a connected component of $G$ and we supposed $G$ has no isolated vertices, we have $G\simeq K_{d,k+1-d}$. Thus $ e(G)=d(k+1-d)\leq{k\choose 2}$.

Otherwise, $G$ is $d:=({k+1\over 2})$-regular. Fix $x\in V$ and set $G':=G-(N(x)\cup\{x\})$. Thus, it is clear that $\ex(G',P_3)+d\leq\ex(G,P_3)<k$, so $\ex(G',P_3)<k-d={k-1\over 2}$. Setting $k':={k-1\over 2}$, we have that $e(G')\leq{k'\choose 2}+\mathbf{1}_{k'=3}$ by induction. Furthermore, as $G$ is triangle-free, $N(x)$ spans no edges, so 

\[
e(G) = e(G')+d^2 \leq \binom{k'}{2} + \mathbf{1}_{k'=3} +d^2 = \frac{1}{8}(3k^2+5) +\mathbf{1}_{k=7} \leq \frac{1}{8}(4k^2-4k) = \binom{k}{2}.
\]

Here, the last inequality follows from $k \geq 4$ necessarily being odd.

\end{proof}

\begin{remark}
For $k \geq 5$, the only extremal graphs in Theorem \ref{thm:p3} are $K_{k-1}^*(1,\dots,1)$ in addition to $K_k$ (when $3 \nmid k$). 
\end{remark}

This can be shown with a little more work with the inductive proof given here. However, the sporadic appearance of $C_4$ and $K_{2,3}$ as extremal examples for $k=3,4$, respectively, overload the case work to establish a base case for the induction, for which we employ NAUTY to do an exhaustive search. The details are excluded here, but can be found at \url{https://github.com/cocox-math/inverse-turan} for the curious reader.

As mentioned before, all $P_3$-free graphs can be built from stars and triangles. So Theorem ~\ref{thm:p3} could be phrased as ``pendant graphs have the smallest (maximum) triangle/star-packing''. This leads naturally to asking which graphs have the smallest star-packing, or in our terminology, determining $\rk_{\{P_3,K_3\}}(k)$.

\begin{theorem}\label{thm:p3k3}
For $\fh = \{ P_3, K_3 \}$, and $k\geq 3$,
\[
\rk_{k}(\fh)=
\begin{cases}
\binom{k+1}{2}-\frac{k+2}{2} & \text{ if $k$ is even;}\\
\binom{k+1}{2}-\frac{k+1}{2} & \text{ if $k$ is odd.}
\end{cases}
\]
Moreover, the only extremal graph for $\rk_{\fh}(k)$ is
\[ G_k:=
\begin{cases}
K_{k+1} \backslash (\frac{k-2}{2}
K_2 \cup P_2)
 & \text{ if $k$ is even;}\\
K_{k+1} \backslash (\frac{k+1}{2}
K_2) & \text{ if $k$ is odd.}
\end{cases}
\]
\end{theorem}

Note that the first of these results has been previously noticed by Ferneyhough, Haas, Hanson and MacGillivray in ~{\cite[Corollary 2]{FHHM02}} using a bound on the domination number due to Vizing ~\cite{V65}. They also used the graphs $G_k$ to provide the lower bounds.
Our new contribution is to show the uniqueness of the extremal graphs. However, said proof is inductive and elementary, so we omit the proof for brevity. The curious reader can find the details at \url{https://github.com/cocox-math/inverse-turan}.

The last small graph we will consider in this paper is $P_1\cup P_2$. Determining $\rk_{P_1\cup P_2}(k)$ will also allow us to completely classify those families of graphs with $\rk_\fh(k)=\infty$, which we will do at the end of this section. As above, it will be important to have a complete classification of $(P_1\cup P_2)$-free graphs.

\begin{lemma}\label{lem:starmatch}
A graph $H$ is $(P_1 \cup P_2)$-free if and only if one of the following holds:
\begin{itemize}
\item $H\simeq sK_2$ for some $s$,
\item $H\simeq K_{1,s}$ for some $s$,
\item $H \subseteq K_4$.
\end{itemize}
\end{lemma}

\begin{proof}
Let $F$ be the line graph of $H$ (whereby $V(F):=E(H)$ and $e_1 \sim_F e_2$ if and only if $e_1$ and $e_2$ share a vertex).
As $H$ is $(P_1 \cup P_2)$-free, for 3 distinct edges $e_1,e_2,e_3 \in E(H)$, then if $e_1\nsim_F e_2$ and $e_2\nsim_F e_3$, then it must be the case that $e_1 \nsim_F e_3$. In particular, the relation $\{(x,y)\in V(F)^2:\text{$x=y$ or $x\nsim_F y$}\}$ is an equivalence relation on $V(F)$, so we may color $V(F)=E(H)$ so that any color class is a matching, and any two edges of a distinct color are incident.
\begin{itemize}
\item Suppose some color class has $s\geq 3$ edges. Since these $s$ edges are disjoint, no other edge can be simultaneously incident to all of these, so every other color class must be empty. Thus $H\simeq sK_2$.
\item Suppose some color class has $2$ edges. Then all other edges must be incident to both of these, so $H \subseteq K_4$.
\item Otherwise, there is $1$ edge in each color, and they are all pairwise incident, so $H\simeq K_3$ or $H\simeq K_{1,s}$ for some $s$.
\end{itemize}
Conversely, all of these graphs are clearly $(P_1 \cup P_2)$-free.
\end{proof}

With this classification, determining $\ex(G,P_1\cup P_2)$ for any graph $G$ is straightforward.

\begin{cor}\label{cor:p2p3degmatch}
For any $G$,
$\ex(G,P_1\cup P_2)=  \max \{\Delta(G),M(G)\}=:t $,
provided $t \geq 6$.
Here $M(G)$ is the size of a maximum matching in $G$.
\end{cor}

\begin{proof}
As any star in $G$ is $(P_1\cup P_2)$-free, certainly $\ex(G,P_1\cup P_2)\geq\Delta(G)$. Similarly, $\ex(G,P_1\cup P_2)\geq M(G)$ as any matching in $G$ is also $(P_1\cup P_2)$-free.

Conversely, take any subgraph $H\subseteq G$ with $t+1>6$ edges, so $H\not\subseteq K_4$. By the definition of $t$, $H$ is neither a star nor a matching, so by Lemma~\ref{lem:starmatch}, $H$ must contain a copy of $P_1\cup P_2$. Therefore, $\ex(G,P_1\cup P_2)\leq t$.
\end{proof}

Using the above Corollary, we can provide lower bounds on $\rk_{P_1\cup P_2}(k)$.

\begin{cor}\label{cor:p2p3}
If $k \geq 7$, then $\rk_{P_1\cup P_2}(k)
\geq \begin{cases}
k^2 - \frac{3}{2} k  & \text{if $k$ is even;}\\
k^2 - k  & \text{if $k$ is odd}.
\end{cases}$
\end{cor}

\begin{proof}
First suppose $k$ is odd and consider $G:=2K_k$, so $\Delta(G)=M(G)=k-1$. Therefore $\ex(G,P_1\cup P_2)<k$ by Corollary~\ref{cor:p2p3degmatch} as $k\geq 7$, so $\rk_{P_1\cup P_2}(k)\geq e(G)=2{k\choose 2}=k^2-k$.

On the other hand, if $k$ is even, start with the Cayley graph $H:=\text{Cay}\bigl(\Z_{2k-1},\bigl[-{k-2\over 2},{k-2\over 2}\bigr]\setminus\{0\}\bigr)$; that is $V(H)=\Z_{2k-1}$ and $xy\in E(H)$ if and only if $x-y\pmod{2k-1}\in\bigl[-{k-2\over 2},{k-2\over 2}\bigr]\setminus\{0\}$. Now, look at all pairs of the form $\{xy:|x-y|={k/2}\}$. Since $k/2$ and $2k-1$ are coprime, these pairs form a Hamilton cycle in the complete graph on $\Z_{2k-1}$, so take any matching $M$ among them of size $k-1$. Finally, consider the graph $G:=(\Z_{2k-1},E(H)\cup M)$. As $M$ and $E(H)$ are disjoint, every vertex of $G$ has degree $(k-1)$ with the exception of one vertex, which has degree $k-2$. Also, $M(G)=k-1$, so $\ex(G,P_1\cup P_2)<k$ again by Corollary~\ref{cor:p2p3degmatch}. Therefore,
\[
\rk_{P_1\cup P_2}(k)\geq e(G)={1\over 2}\sum_{v\in V(G)}\deg(v)={1\over 2}\bigl((2k-2)(k-1)+(k-2)\bigr)=k^2-{3\over 2}k.\qedhere
\]
\end{proof}

To yield upper bounds on $\rk_{P_1\cup P_2}(k)$, we prove a general bound on the number of edges of a graph based on its maximum degree and matching number. A similar theorem was proved by Abbot, Hanson and Sauer~\cite{AHS72} in the context of the Erd\H{o}s--Rado sunflower lemma, but we provide a full proof for completeness.

\begin{theorem}\label{thm:degmatch}
For a graph $G$, $ e(G)\leq(\Delta(G)+1)M(G)$.
Furthermore, the inequalities in Corollary~\ref{cor:p2p3} are in fact equalities.
\end{theorem}

In order to prove this, we will need the following proposition, which is an immediate consequence of the Gallai--Edmonds decomposition of a graph (c.f.\ \cite{LP09} pp.\ 93--95).

\begin{prop}\label{prop:gallai}
If $G$ is a connected graph with the property that for every $v\in V$, $M(G-v)=M(G)$, then $G$ has an odd number of vertices and $M(G)={|V(G)|-1\over 2}$.
\end{prop}

\begin{proof}[Proof of Theorem~\ref{thm:degmatch}]
Let $G$ be a graph with $M(G)\leq M$ and $\Delta(G)\leq\Delta$. Suppose $G$ has components $S_1,\dots,S_s,H_1,\dots,H_t$, where $S_i$ is a star of degree at most $\Delta$. We will consider a series of reductions of $G$ that maintain the matching and degree restrictions and not decrease the number of edges. We first claim that for each $i$ and any $v\in V(H_i)$, we may suppose that $M(H_i-v)=M(H_i)$. To see this, suppose that this is not the case for some $i$ and $v$. In this case, let $G'$ be the graph formed by replacing $H_i$ with $H_i'=H_i-v$ and adding a copy of $K_{1,\Delta}$. As $\deg(v)\leq\Delta$, we have $\Delta(G') = \Delta$ and $e(G')\geq e(G)$. Furthermore, since every maximum matching in $H_i$ used $v$, $M(H_i')=M(H_i)-1$, so as $M(K_{1,\Delta})=1$, we have $M(G')=M(G)\leq M$. Thus, we may assume that $M(H_i-v)=M(H_i)$ for all $i$ and $v\in V(H_i)$. 

As such, $|V(H_i)|$ is odd and $M(H_i)={|V(H_i)|-1\over 2}$ by Proposition~\ref{prop:gallai}. We now claim that we may suppose that $|V(H_i)|\geq\Delta+1$ for all $i$. If not, form $G'$ by replacing $H_i$ with a copy of ${|V(H_i)|-1\over 2}K_{1,\Delta}$. Clearly $\Delta(G') = \Delta$, and $M(G')=M(G)$ by the previous comment. Finally, 
\[
e(G')- e(G)={|V(H_i)|-1\over 2}\Delta-e(H_i)\geq{|V(H_i)|-1\over 2}|V(H_i)|-{|V(H_i)|\choose 2}=0,
\]
so we may suppose this property of $G$. Additionally, as $|V(H_i)|$ is odd, this property tells us $|V(H_i)|\geq\Delta+1+\mathbf{1}_{\text{$\Delta$ odd}}$.

Now, 
\[
M\geq M(G)=s+{1\over 2}\sum_{i=1}^t\bigl(|V(H_i)|-1\bigr),
\]
so we find 
\[
t\leq \bigg\lfloor{2M\over \min_i\{|V(H_i)|-1\}}\bigg\rfloor\leq \bigg\lfloor{2M\over \Delta+\mathbf{1}_{\text{$\Delta$ odd}}}\bigg\rfloor.
\]
 Rewriting the above equation as $s+{1\over 2}\sum_{i=1}^t|V(H_i)|\leq M+t/2$, we calculate
\begin{align}\label{eqn:improve}
 e(G) &= \sum_{i=1}^s e(S_i)+\sum_{i=1}^te(H_i) \nonumber\\
&\leq s\Delta+{\Delta\over 2}\sum_{i=1}^t|V(H_i)|\nonumber\\
&\leq \Delta\biggl(M+{t\over 2}\biggr)\nonumber\\
& \leq \Delta\biggl(M+{1\over 2}\bigg\lfloor{2M\over\Delta+\mathbf{1}_{\text{$\Delta$ odd}}}\bigg\rfloor\biggr)  \\
&\leq \bigl(\Delta+1\bigr)M.\nonumber
\end{align}
Now, take any $k \geq 7$ and let $G$ be a graph with $\ex(G,P_1\cup P_2)<k$, so we must have $\Delta(G), M(G) \leq k-1$. 
Now, when $k$ is odd, immediately $ e(G) \leq (\Delta(G)+1)M(G) \leq k(k-1)$; hence $\rk_{P_1\cup P_2}(k)=k^2 - k$.\\
 When $k$ is even, note that either $\Delta \leq k-2,$ in which case immediately
 $ e(G) \leq (k-1)^2 < k^2 - \frac{3}{2}k,$
 or else $\Delta =k-1$, so by Equation~\eqref{eqn:improve},
 \[
  e(G)\leq(k-1)\biggl((k-1)+{1\over 2}\bigg\lfloor {2(k-1)\over k}\bigg\rfloor\biggr)=(k-1)\biggl(k-{1\over 2}\biggr)=k^2-{3\over 2}k+{1\over 2}.
 \]
 Thus as $k$ is even, we have $ e(G)\leq k^2-{3\over 2}k$, so $\rk_{P_1\cup P_2}(k)=k^2-{3\over 2}k$ in this case.
\end{proof}

We conclude this section with an absolute upper bound on the finite values of $\rk_\fh(k)$. 

\begin{cor}\label{cor:2bound}
If $\fh$ is a family of simple 2-uniform graphs which contains neither a star nor a matching, then $\rk_\fh(k)\leq\rk_{P_1\cup P_2}(k)\leq k(k-1)$ for $k\geq 3$.
\end{cor}

\begin{proof}
	If $H$ is a simple 2-uniform graph, then either $H$ is a star, a matching, or contains one of $K_3,P_3,P_1\cup P_2$. As such,
	\[
	\rk_\fh(k)\leq\max\bigl\{\rk_{K_3}(k),\rk_{P_3}(k),\rk_{P_1\cup P_2}(k)\bigr\}=\rk_{P_1\cup P_2}(k)\leq k(k-1),
	\]
	which follows from the fact that $\rk_{K_3}(k)\leq 2k$ (Theorem~\ref{thm:cliques}), $\rk_{P_3}(k)={k\choose 2}+\mathbf{1}_{k=3}$ (Theorem~\ref{thm:p3}) and $\rk_{P_1\cup P_2}(k)= k(k-1)$ if $k$ is odd and $\rk_{P_1\cup P_2}(k)=k\bigr(k-{3\over 2}\bigr)$ if $k$ is even (Theorem~\ref{thm:degmatch}).
\end{proof}

\subsection{The Trichotomy of the Extremal Numbers}\label{sec:tri}

We end our discussion of 2-uniform graphs by showing how the same ideas as in Theorem~\ref{thm:eshypergraphs} and~\ref{thm:cliques} can be used to give non-trivial bounds on $\rk_H(k),\rk_H^*(k)$ when $H$ is not only non-bipartite. Note that a similar statement holds for hypergraphs of higher uniformity.

\begin{theorem}\label{thm:2degen}
	Let $H$ be a simple $2$-uniform graph which is neither a matching nor a star.
	\begin{enumerate}
		\item If $\ex(K_n,H)\leq O(n^\alpha)$, then $		\rk_H(k)\geq\Omega(k^{2/\alpha})$, and
		\item If $\ex(K_n,H)\geq\Omega(n^\beta)$, then $
		\rk_H^*(k)\leq O(k^{3-\beta})$.
	\end{enumerate}

\end{theorem}
\begin{proof}
	This proof is very similar to that of Theorem~\ref{thm:eshypergraphs}.
	\begin{enumerate}
		\item Set $k=\ex(K_n,H)+1\leq O(n^\alpha)$, so $K_n$ has ${n\choose 2}\geq \Omega(k^{2/\alpha})$ edges and $\ex(K_n,H)<k$, so $\rk_H(k)\geq\Omega(k^{2/\alpha})$.
		
		\item If $H$ is a forest which is neither a matching nor a star, then $\ex(K_n,H)=\Theta(n)$ and $\rk_H^*(k)\leq O(k^2)$ by Corollary~\ref{cor:2bound}. Thus, suppose $H$ has a cycle.
		
		First, suppose $H$ is connected and let $G$ be an extremal graph for $\rk_H^*(k)$. Since $H$ is connected, we may suppose that $G$ is also connected by Corollary~\ref{cor:connect}. Therefore, since $\ex(G,H)<k$ and $H$ contains a cycle, we must have $|V(G)|\leq k$ since any spanning tree of $G$ is $H$-free. Set $n=|V(G)|$ and let $F\subseteq K_n$ be an extremal graph for $\ex(K_n,H)$. Uniformly selecting a random copy $F'$ of $F$ from $K_n$ and setting $F^*=\{e\in E(G):e\in E(F')\}$, where any multiedges are preserved, we see that $F^*$ is $H$-free, so
		\[
		k>\E e(F^*)={\ex(K_n,H)\over{n\choose 2}}\cdot e(G)\geq \Omega\bigl(e(G)n^{\beta-2}\bigr).
		\]
		Thus, $\rk_H^*(k)=e(G)\leq O(k^{3-\beta})$.
		
		Now, if $H$ has a cycle and is not connected, let $H_1,\dots,H_\ell$ be the connected components of $H$. Without loss, suppose that $\ex(K_n,H_1)=\max_i\ex(K_n,H_i)$, and notice that $\ex(K_n,H_1)=\bigl(1+o(1)\bigr)\ex(K_n,H)\geq\Omega(n^\beta)$ since $H_1$ contains a cycle. Finally, since $H_1\subseteq H$, we conclude that $\rk_H^*(k)\leq\rk_{H_1}^*(k)\leq O(k^{3-\beta})$.\qedhere 
	\end{enumerate}
\end{proof}

From this, we recover a similar trichotomy to that of the Tur\'an number. Recall that if $H$ is a 2-uniform simple graph with at least two edges, then
\begin{enumerate}
	\item $\ex(K_n,H)=\Theta(n^2)$ if $\chi(H)\geq 3$,
	\item $\Omega(n^\alpha)\leq\ex(K_n,H)\leq O(n^\beta)$ for some $1<\alpha\leq\beta<2$ if $\chi(H)=2$ and $H$ has a cycle, and
	\item $\ex(K_n,H)=\Theta(n)$ if $H$ is a forest.
\end{enumerate}
In the case of $\rk_H(k)$, we have:
\begin{cor}
	Let $H$ be a simple 2-uniform graph.
	\begin{enumerate}
		\item $\rk_H(k),\rk_H^*(k)=\Theta(k)$ if $\chi(H)\geq 3$,
		\item $\Omega(k^\alpha)\leq\rk_H(k)\leq\rk_H^*(k)\leq O(k^\beta)$ for some $1<\alpha\leq\beta<2$ if $\chi(H)=2$ and $H$ has a cycle, and
		\item $\rk_H(k),\rk_H^*(k)=\Theta(k^2)$ if $H$ is a forest which is neither a star nor a matching.
	\end{enumerate}
\end{cor}


Going back to Theorem~\ref{thm:2degen}, since $\ex(K_n,C_4)=\Theta(n^{3/2})$, this immediately implies the following:
\begin{cor}
	$\Omega(k^{4/3})\leq\rk_{C_4}(k)\leq\rk_{C_4}^*(k)\leq O(k^{3/2})$.
\end{cor}
Note that the upper bound of $k^{3/2}$ can also be attained via Proposition~\ref{prop:genupper}. In fact, $\ex(K_{q^2+q+1}, C_4)=\frac{1}{2}q(q+1)^2$
and $\ex(K_{q^2+q+1,q^2+q+1},C_4)=(q^2+q+1)(q+1)$
~\cite{FZ96} when $q$ is a prime power.
So setting $k=\frac{1}{2}q(q+1)^2+1$ or $k=(q^2+q+1)(q+1)+1$ shows $\rk_{C_4}(k) \gtrsim 2^{1/3}k^{4/3} $ and $\rk_{C_4}(k) \gtrsim k^{4/3}$ for these respective values of $k$,
 suggesting cliques even outperform complete bipartite graphs when forcing $C_4$'s.
This leads us to the following general question.
\begin{question}
	Let $H$ simple 2-uniform graph which is neither a star nor a matching. If $\ex(K_n,H)=\Theta(n^\beta)$, then must it be the case that $\rk_H(k),\rk_H^*(k)=\Theta(n^{2/\beta})$?
\end{question}
As we have shown, the only possible graphs $H$ for which the answer to the above question can be negative are bipartite graphs which contain a cycle. However, for these $H$'s, we do not even have a good understanding of $\ex(K_n,H)$ to begin with.

\section{Non-uniform Hypergraphs}\label{sec:nonu}

In this section, we will consider an analogue of $\rk_\fh(k)$ and $\rk_\fh^*(k)$ when the graphs in $\fh$ are not uniform. We say that a (multi)hypergraph $H$ is a non-uniform graph if it has two edges of different uniformities. $\fh$ is then said to be a family of non-uniform graphs if each $H\in\fh$ is non-uniform. Notice that $\rk_\fh(k)$ and $\rk_\fh^*(k)$ are actually already well-defined if $\fh$ is a family of non-uniform graphs; namely, for example, $\rk_\fh^*(k)$ is still the maximum number of edges in a (multi)graph $G$ such that any $k$ edges of $G$ contains a copy of some $H\in\fh$. The only possible difference is that here $G$ must be non-uniform itself. The distinction between simple non-uniform graphs and non-uniform (multi)graphs is simply that multigraphs are allowed to have parallel edges. Note that $G$ can have edges $e,f$ with $e\subsetneq f$ and still be considered simple non-uniform graph provided it does not have any repeated edges. In this section, a ``graph'' is allowed to be a non-uniform multigraph.

Throughout this section, for a graph $G$, we will use $E_i(G):=\{e\in E(G):|e|=i\}$, that is $E_i(G)$ is the set of edges of uniformity $i$, and write $e_i(G)=|E_i(G)|$.

\begin{prop}
If $H$ is a non-uniform (multi)graph, then $\rk_H^*(k)\leq 2(k-1)$. Additionally, if $\fH$ is a finite family of non-uniform graphs, then $\rk_\fh^*(k)$ is always finite.
\end{prop}

\begin{proof}
Since $H$ is non-uniform, there is some $r\neq s$ with $E_r(H),E_s(H)\neq\varnothing$. Now, let $G$ be any graph with $\ex(G,H)<k$. As any $F\subseteq G$ with $E_r(F)=\varnothing$ or $E_s(F)=\varnothing$ is trivially $H$-free, we know that $e_r(G)<k$ and $e(G)-e_r(G)<k$, therefore, $e(G)\leq 2(k-1)$.

Similarly, for a finite family of non-uniform graphs $\fH$, let $U=\{i\in\Z:E_i(H)\neq\varnothing\text{ for some $H\in\fH$}\}$. Let $G$ be a graph with $\ex(G,\fH)<k$; certainly we may suppose that the edges in $G$ are only of the sizes in $U$. Thus, by the same argument as above, $e_i(G)<k$ for all $i\in U$ since each $H\in\fH$ is non-uniform, so $e(G)\leq |U|(k-1)$, which is finite since $\fH$ consisted only of finitely many graphs.
\end{proof}

We quickly remark that $\rk_\fh^*(k),\rk_\fh(k)$ are not necessarily finite when $\fH$ is not of finite size. Namely, for positive integers $r,t$, let $H_{r,t}$ be the non-uniform graph consisting of two disjoint edges $e,s$ with $|e|=r$, $|s|=t$. Then $\fH=\{H_{r,t}:1\leq r<t\}$ has $\rk_\fh(k)=\infty$ when $k\geq 2$, as is realized by taking a host graph with disjoint edges $e_1,\dots,e_s$ where $|e_i|=i$.

We now turn our attention to a non-uniform graph which yields a surprising answer to $\rk_H^*(k)$; namely $\rk_H^*(k)\sim xk$ where $x$ is an irrational number.

The dumbbell graph, denoted $\cal{D}$, is the non-uniform graph consisting of a single 2-uniform edge with one loop at each vertex.

\begin{theorem}\label{thm:necklace}
	$\rk_\cal{D}^*(k)=\bigl(\phi-o(1)\bigr)k$ where $\phi=1.618\ldots$ is the golden ratio.
\end{theorem}

\begin{proof}
To simplify notation, we set $\varphi=1/\phi$. Note that $\varphi$ satisfies $\varphi^2+\varphi=1$.

\textit{Upper bound.} Let $G$ be any (multi)graph with $\ex(G,\cal{D})<k$; certainly we may assume that $G$ contains only edges of uniformities $1$ and and $2$, that is only has 2-uniform edges and loops. Let $V'\subseteq V(G)$ be formed by including each vertex of $V(G)$ into $V'$ with probability $\varphi$, and form $G'\subseteq G$ by taking any loops on a vertex of $V'$ along with any 2-uniform edge which is not completely contained in $V'$. By construction, $G'$ is $\cal{D}$-free, so
\[
k>\E e(G')=\varphi\cdot e_1(G)+(1-\varphi^2)e_2(G)=\varphi\bigl(e_1(G)+e_2(G)\bigr)=\varphi\cdot e(G),
\]
so $e(G)<k/\varphi=\phi\cdot k$.

\textit{Lower bound.} We will show $\rk_{\cal{D}}^*(k) \geq \phi\cdot k - O(k^{2/3})$.

Fix a large $k$, and construct the multigraph $G$ on $n=\Theta(k^{1/3})$ vertices with:

\begin{itemize}
\item $\displaystyle t:= \biggl\lfloor {k\over{n\choose 2}}\cdot {1\over \varphi+2\varphi^2} \biggr\rfloor$ parallel 2-uniform edges spanning every pair of vertices, and
 \item

$ \displaystyle s:= \biggl\lfloor {k\over n} \cdot {2\varphi-{1\over n}\over \varphi+2\varphi^2} \biggr\rfloor
$ 
 loops at each vertex.\footnote{The constants $1/(\varphi+2\varphi^2)$ and $2\varphi/(\varphi+2\varphi^2)$ may be found by solving the natural linear program, but this derivation is not necessary for the proof.}
\end{itemize}

This way, $G$ has
\begin{align*}
	e(G) &= t{n\choose 2}+sn \geq \biggl({k\over {n\choose 2}}\cdot {1\over\varphi+2\varphi^2}-1\biggr){n\choose 2}+\biggl({k\over n}\cdot{2\varphi-{1\over n}\over \varphi+2\varphi^2}-1\biggr)n\\
	&={k\over\varphi}-{n\choose 2}-n-{1\over\varphi+2\varphi^2}{k\over n}=\phi\cdot k-O(k^{2/3}).
\end{align*}

We now wish to show that $\ex(G,\cal{D})<k$, so take any $\cal{D}$-free subgraph $H \subseteq G$; we will show $e(H) <k$.

Let $L\subseteq V(H)$ be the vertices of $H$ with at least one loop. Write $\beta n:= |L|$, then certainly $H$ has 
\[
e_1(H)\leq\beta n s\leq {k\over \varphi+2\varphi^2}\biggl(2\phi\beta-{\beta\over n}\biggr).
\]
Furthermore, since $H$ is $\cal{D}$-free, $H$ cannot have a 2-edge completely contained in $L$, so 
\[
e_2(H)\leq t\biggl({n\choose 2}-{\beta n\choose 2}\biggr)\leq {k\over \varphi+2\varphi^2}\biggl(1-\beta^2+{\beta\over n}\biggr).
\]

We also note that $\beta\neq\varphi$ since $\beta$ is rational and $\varphi$ is irrational; therefore, by the mean value theorem, there is some $\theta$ strictly between $\varphi$ and $\beta$ such that $2(\varphi-\beta)\theta=\varphi^2-\beta^2$. Hence,
\[
\varphi+\varphi^2+2\beta\varphi-\beta^2=\varphi+2\varphi^2+2(\varphi-\beta)(\theta-\varphi)<\varphi+2\varphi^2,
\]
since either $\varphi<\theta<\beta$ or $\beta<\theta<\varphi$. Putting this together,
\begin{align*}
	{e(H)\over k} &= {e_1(H)+e_2(H)\over k}\leq{1\over \varphi+2\varphi^2}\biggl(2\beta\varphi-{\beta\over n}+1-\beta^2+{\beta\over n}\biggr)\\
	&={1+2\beta\varphi-\beta^2\over \varphi+2\varphi^2}={\varphi+\varphi^2+2\beta\varphi-\beta^2\over \varphi+2\varphi^2}<{\varphi+2\varphi^2\over\varphi+2\varphi^2}=1.\qedhere
\end{align*}
\end{proof}

\begin{remark}
	Let $\cal{D}_r$ denote the non-uniform graph consisting of a single $r$-uniform edge with one loop at each vertex. The above proof can be generalized to show that $\rk_{\cal{D}_r}^*(k)\sim \phi_r\cdot k$ where $\phi_r$ is the unique positive solution to $X^r-X^{r-1}=1$.
\end{remark}

In contrast to Conjecture~\ref{conj:multi}, we now prove that $\rk_H(k)$ and $\rk_H^*(k)$ can differ if $H$ is a non-uniform graph.

\begin{theorem}\label{thm:multibeats}
If $G$ is a graph with $1$-uniform edges and $2$-uniform edges, where each vertex has at most one loop (but any $2$-uniform edges can have higher multiplicity), then $\ex(G,\cal{D})\geq{2\over 3}e(G)$.
\end{theorem}
\begin{proof}
Let $G$ be a graph with only $2$-uniform edges and loops where each vertex has at most one loop. As before, let $E_i(G)$ denote the set of $i$-uniform edges, so $E(G)=E_1(G)\cup E_2(G)$. We begin by claiming that we may suppose that every vertex of $G$ has a loop. If some $v\in V(G)$ does not have a loop, then either $v$ is isolated, in which case we may simply delete $v$, or $v$ is incident to some $e\in E_2(G)$. Let $G'$ be formed by deleting $e$ and adding a loop around $v$. Certainly $\ex(G',\cal{D})\leq\ex(G,\cal{D})$ as the edge $e\in E_2(G)$ cannot be used in any copy of $\O_2$. After this reduction, we know that $ e(G)=e_1(G)+e_2(G)=|V(G)|+e_2(G)$. 
Additionally, we may suppose that every vertex is incident to some $e\in E_2(G)$. To see this, suppose $v\in V(G)$ is not incident to any edge in $E_2(G)$; pick any $e\in E_2(G)$ and form $G'$ by removing the loop from $v$ and adding an additional copy of the edge $e$. As the loop around $v$ cannot be used in any copy of $\cal{D}$ in $G$, we see that $\ex(G',\cal{D})\leq\ex(G,\cal{D})$.

We now prove the statement by induction on the number of odd cycles in $G$ (that is, the number of odd cycles in $E_2(G)$).

For the base case, suppose that $E_2(G)$ is bipartite with partite sets $A,B$ where $|A|\geq|B|$. In this case, if we take every edge in $E_2(G)$ and every loop around a vertex in $A$, we end up with an $\cal{D}$-free graph as no two loops are joined by an edge. Now, as $G$ has no vertices not incident to a $2$-edge, we must have $e_2(G)\geq|B|$, so as $|A|\geq|B|$, we have
\[
\ex(G,\cal{D})\geq e_2(G)+|A|\geq{2\over 3} e(G).
\]

Now suppose that $E_2(G)$ is not bipartite and let $C\subseteq G$ be an induced copy of $C_{2t+1}$ for some $t$, possibly with some multiedges. Set $G':=G\setminus C$.

Now, for a fixed set of vertices $S \subseteq V(C)$, we may form $H_S \subseteq G$ by collecting together the following edges:

\begin{itemize}
\item All 2-edges in the cycle $C$ itself (there are at least $2t+1$ of these),
\item The 2-edges from $C\backslash S$ to $V \backslash C$, 
\item All loops in $S$,
\item $E(H')$ for some extremal $\cal{D}$-free $H' \subseteq G'$.
\end{itemize}
Provided $S$ contains no two adjacent vertices in $C$, $H_S$ is $\cal{D}$-free.

\begin{figure}[ht]
\centering
\begin{tikzpicture}
\draw [thick] (-0.1,0) circle [radius =1.2];
\draw [ultra thick,red]  (0.2,0) node {${\color{black} H'}$} circle [radius =0.8];
\filldraw (-0.1,1.5) node {$G'$};

\foreach \a in {1,2,...,5}
 { \filldraw (3,0)+(\a*72:1) circle [radius =0.045];
 \draw [red, ultra thick] ($(3,0)+(\a*72:1)$)--($(3,0)+(\a*72+72:1)$);
 \draw [thick] ($(3,0)+(\a*72:1.125)$) circle [radius =0.15];
 }
 
\draw [red,ultra thick] ($(3,0)+(72:1)$) to [out=270, in=180] ($(3,0)+(0:1)$);
 
 \draw [red,ultra thick] ($(3,0)+(72:1.125)$) circle [radius =0.15];
  \draw [red,ultra thick] ($(3,0)+(288:1.125)$) circle [radius =0.15];
 
 \draw [thick] ($(3,0)+(72:1)$) --(-0.1,1.2);
 \draw [thick] (-0.1,-1.2)--($(3,0)+(288:1)$);
 
 \draw [red,ultra thick] ($(0.5,0) +(15:1)$)--($(3,0)+(144:1)$)--
 ($(0.5,0) +(45:1.2)$);
  \draw [red,ultra thick] ($(0.5,0) +(345:1)$)--($(3,0)+(216:1)$)--
 ($(0.5,0) +(315:1.2)$);
 
 \foreach \a in {1,2,...,5}
 { \filldraw (3,0)+(\a*72:1) circle [radius =0.08];

 }

\filldraw (3,0) node {$C$}; 
\end{tikzpicture}
\caption{\label{fig:O2} $H_S$ edges in red, $G\backslash H_S$ edges in black (here $|S|=2)$.}
\end{figure}

So, suppose we choose $S \subseteq V(C)$ by picking an independent set of size $\bigl\lceil \frac{2t+1}{3} \bigr\rceil$ uniformly at random with probability $\frac{2t+1}{3}-\bigl\lfloor \frac{2t+1}{3}\bigr\rfloor$, and
otherwise an independent set of size
$\bigl\lfloor \frac{2t+1}{3} \bigr\rfloor$ uniformly at random. This is done so that for every $v\in V(C)$,
\begin{align*}
\Pr[v\in S] &=\biggl({2t+1\over 3}-\biggl\lfloor {2t+1\over 3}\biggr\rfloor\biggr){\bigl\lceil{2t+1\over 3}\bigr\rceil\over 2t+1 }+\biggl(1-{2t+1\over 3}+\biggl\lfloor {2t+1\over 3}\biggr\rfloor\biggr){\bigl\lfloor{2t+1\over 3}\bigr\rfloor\over 2t+1}\\
&= {\bigl\lfloor{2t+1\over 3}\bigr\rfloor\over 2t+1}\biggl(1-\biggl\lceil{2t+1\over 3}\biggr\rceil+\biggl\lfloor{2t+1\over 3}\biggr\rfloor\biggr)+{1\over 3}\biggl(\biggl\lceil{2t+1\over 3}\biggr\rceil-\biggl\lfloor{2t+1\over 3}\biggr\rfloor\biggr)\\
&={1\over 3}.
\end{align*}
In particular, $\E|S|={1\over 3}(2t+1)$. Finally, recalling that by induction, $e(H')=\ex(G',\cal{D})\geq{2\over 3}e(G')$, we have in total
\begin{align*}
\E e(H_S)
&=e_2(C)
 + \frac{2}{3} e[C ,V \backslash C] +\frac{1}{3}(2t+1) +
 \frac{2}{3} e(G')\\
& \geq 
 \frac{2}{3}e_2(C)
 + \frac{2}{3} e[C ,V \backslash C]  + \frac{2}{3}(2t+1)
+ \frac{2}{3} e(G')
 =\frac{2}{3}  e(G).
 \end{align*}
 So some such $S$ yields an $\cal{D}$-free $H_S$ with at least this many edges, as desired.
\end{proof}

Now, notice that the above proves that $\rk_\cal{D}(k)\leq\bigl\lfloor{3\over 2}(k-1)\bigr\rfloor$. On the other hand, let $G$ be the non-uniform graph which consists of $n$ disjoint copies of $\cal{D}$ if $k=2n+1$, or the graph that consists of $n-1$ disjoint copies of $\cal{D}$ along with a single isolated edge if $k=2n$. In either case, $\ex(G,\cal{D})<k$ and $e(G)=\bigl\lfloor{3\over 2}(k-1)\bigr\rfloor$, so we have proved Theorem~\ref{thm:dumbbell}, showing that the analogue of Conjecture~\ref{conj:multi} can fail for non-uniform graphs.

\begin{question}
	What is $\rk_{\cal{D}_r}(k)$ for $r\geq 3$?
\end{question}

%

\section{1-Uniform Graphs}\label{sec:1u}

A $1$-uniform multigraph on $n$ vertices is equivalent to its degree sequence $(d_1,\dots,d_t)$ where $d_i$ is the number of loops at vertex $i$. For 1-uniform graphs $H=(d_1,\dots,d_t)$ and $G=(x_1,\dots,x_n)$, $H\subseteq G$ if and only if there is an injection $f:[t]\to[n]$ such that for every $i\in[t]$, $d_i\leq x_{f(i)}$.

We quickly note that a $1$-uniform graph $H$ is a sunflower if and only if it is of the form $H=(1,1,\dots,1)$ or $H=(r)$ for some $r$. Because of this, looking at $\rk_H(k)$ is uninteresting since every simple 1-uniform graph is a sunflower, so we look only at $\rk_H^*(k)$ when $H=(d_1,\dots,d_t)$ where $d_1\geq\dots\geq d_t\geq 1$ and $d_1,t\geq 2$.

One reason for caring about 1-uniform graphs in this context is that it also settles the question for multi-stars. For positive integers $d_1,\dots,d_t$, the multi-star $S_{d_1,\dots,d_t}$ is a star on $t+1$ vertices whose edges have multiplicities $d_1,\dots,d_t$.

\begin{lemma}
	For positive integers $d_1,\dots,d_t$, if $H=(d_1,\dots,d_t)$, then $\rk_{S_{d_1,\dots,d_t}}^*(k)=\rk_H^*(k)$.
\end{lemma}

\begin{proof}
	Since every edge of a multi-star is uniquely determined by its leaf-vertex, we observe that
	\[ \ex(S_{x_1,\dots,x_n},S_{d_1,\dots,d_t})=\ex\bigr((x_1,\dots,x_n),H\bigr);
	\]
	hence, $\rk_{S_{d_1,\dots,d_t}}^*(k)\geq\rk_H^*(k)$.
	
	On the other hand, let $G$ be a 2-uniform (multi)graph with $\ex(G,S_{d_1,\dots,d_t})<k$ where $G$ has edges with multiplicites $x_1,\dots,x_n$. Notice that $\ex(G,S_{d_1,\dots,d_t})\geq\ex(S_{x_1,\dots,x_n},S_{d_1,\dots,d_t})$ since any copies of $S_{d_1,\dots,d_t}$ are preserved upon enforcing all non-parallel edges to share a common vertex. Thus, 
	\[
	\ex\bigl((x_1,\dots,x_n),H\bigr)=\ex(S_{x_1,\dots,x_n},S_{d_1,\dots,d_t})<k,
	\]
	and $e(S_{x_1,\dots,x_n})=e(G)$, so we have $\rk_H^*(k)\geq\rk_{S_{d_1,\dots,d_t}}^*(k)$.
\end{proof}

\begin{theorem}\label{thm:1u}
For every $1$-uniform multigraph $H=(d_1,d_2,\dots,d_t)$ with $d_1\geq\dots\geq d_t\geq 1$ where $d_1,t\geq 2$, there exists a constant $c_H$ such that $\rk_H^*(k)=\bigl(c_H+o(1)\bigr)k^2$. Additionally, $c_H$ can be determined in polynomial time and satisfies ${1\over 4(t-1)(d_1-1)}\leq c_H\leq{1\over (t-1)(d_1-1)}$.
\end{theorem}

\begin{proof}
Let $H=(d_1,\dots,d_t)$ where $d_1\geq\dots\geq d_t \geq 1$ and $d_1,t\geq 2$.
%
%
%
%

We note that $F\subseteq G$ with $F=(f_1,\dots,f_n)$ can be assumed to have $f_1\geq\dots\geq f_n$. Thus, it is clear that $F$ is $H$-free if and only if there is some $t'\in[t]$ such that $f_{t'}<d_{t'}$. Thus, for $t'\in[t]$, let $G_{t'}=(x_1',\dots,x_n')$ where $x_i'=x_i$ for $i<t'$ and $x_i'=\min\{x_i,d_{t'}-1\}$ for all $i\geq t'$. By the earlier note, $G_{t'}$ is $H$-free for every $t'\in[t]$; furthermore, any $F\subseteq G$ that is $H$-free must be contained in some $G_{t'}$. Thus,
\[
\ex(G,H)=\max_{t'\in[t]}e(G_{t'}).
\]
Since $H\neq(1,1,\dots,1)$, we know that if $\ex(G,H)<k$, then $n\leq k-1$. Thus, we may formulate the following non-linear integer program for $\rk_H^*(k)$:
\[
\begin{array}{cccll}
\rk_H^*(k) & = & \max & \sum_{i=1}^{k-1} x_i &\\
&& \text{s.t.} & \sum_{i=1}^{t'-1}x_i+\sum_{i=t'}^{k-1}\min\{x_i,d_{t'}-1\}\leq k-1 & \text{for all $t'\in[t]$}\\
&&& x_i\in \Z_{\geq 0} & \text{for all $i\in[k-1]$}.
\end{array}
\]


Fix a feasible $G$.
Note that, since $t\geq 2$, taking $t'=2$ shows $x_1 \leq x_1 + \sum_{i=2}^{k-1} \min\{x_i,d_2-1\} \leq k-1$. 

Now, define $j:= \max \{i: x_i \geq d_1 \}$.
Then $\sum_{i>j} x_i \leq (k-1)(d_1-1)\leq d_1k$.
Furthermore,
if the largest $j$ vertices $(x_1, \dots, x_j)$ differ in degree by $\geq 2$, then certainly $x_i \geq x_{i+1}+1 \geq\dots\geq x_{\ell-1}+1 \geq  x_\ell +2$ for some $i < \ell \leq j$.

Then $G'$ formed by replacing $x_i,x_\ell$ with $x_i-1,x_\ell + 1$ respectively (noting the degree sequence is still decreasing) is still feasible, for otherwise
the first condition is violated for some $t'\in[t]$. This would mean
\[
\min\{ x_i -1,d_{t'}-1 \} + \min \{ x_\ell +1, d_{t'} -1 \} > \min\{ x_i ,d_{t'}-1 \} + \min \{ x_\ell , d_{t'} -1 \},\]
as only $x_i$ and $x_\ell$ changed in value when forming $G'$.
 Thus
$\min \{ x_{\ell} +1,d_{t'}-1\} =x_{\ell}+1$ so $x_{\ell} \leq d_{t'}-2 <d_1 $; a contradiction.

Thus, we may suppose $G=(x_1,\dots,x_n)$ where $|x_i-x_\ell|\leq 1$ for all $i,\ell\leq j$. From this, define $G^{(1)}:=(\underbrace{x_j,x_j,\dots,x_j}_j,0,\dots,0)$, which is also feasible and has
\[
 e(G)-e(G^{(1)})=\sum_{i=1}^j(x_i-x_j)+\sum_{i=j+1}^{k-1}x_i\leq j+d_1k=O(k).
\]
As such, we have $f^{(1)}(H)\leq\rk_H^*(k)\leq f^{(1)}(H)+O(k)$ where
%
\[
\begin{array}{cccll}
f_k^{(1)}(H) & = & \max & jx &\\
&& \text{s.t.} & (t'-1)x+\sum_{i=t'}^j\min\{x,d_{t'}-1\}\leq k-1 & \text{for all $t'\in[t]$}\\
&&& x,j\in\Z_{\geq 0},\ x \geq d_1,
\end{array}
\]
where the lower bound follows from the fact that for a feasible pair $(x,j)$, the 1-graph $(\underbrace{x,x,\dots, x}_j, 0, \dots, 0)$ satisfies the original program.

To simplify further, note that $x> d_{t'}-1$ for any $t'\in[t]$ as $x\geq d_1$, so we know that $\min\{x,d_{t'}-1\}=d_{t'}-1$. Furthermore, if a feasible $(x,j)$ has $j<t,$  then the objective is $xj<(k-1)t=O(k)$.
Whether or not the optimum is among such $(x,j)$, this shows we decrease the objective by at most $O(k)$
 upon imposing the restriction $j \geq t$. Thus,
\[
\begin{array}{cccll}
f_k^{(2)}(H) & = & \max & jx &\\
&& \text{s.t.} & (t'-1)x+(j-t'+1)(d_{t'}-1)\leq k-1 & \text{for all $t'\in[t]$}\\
&&& x \geq d_1\\
&&& j \geq t\\
&&& x,j\in\Z,
\end{array}
\]
has $f_k^{(2)}(H)\leq f_k^{(1)}(H)\leq f_k^{(2)}(H)+O(k)$.

Next, replace $j$ with $j-t$, and $x$ with $x-d_1$, noting that the objective function decreases by $xj - (x-d_1)(j-t)\leq xt+jd_1\leq O(k)$. Thus
\[
\begin{array}{cccll}
f_{k}^{(3)}(H) & = & \max & jx &\\
&& \text{s.t.} & (t'-1)x+(d_{t'}-1)j\leq k-1 & \text{for all $t'\in[t]$}\\
&&& x,j\in\Z_{\geq 0}
\end{array}
\]
satisfies $f_k^{(3)}(H)\leq f_k^{(2)}(H)\leq f_k^{(3)}(H)+O(k)$.

We now relax the integrality of $x,j$ to attain
\[
\begin{array}{cccll}
f^{(4)}_k(H) & = & \max & jx &\\
&& \text{s.t.} & (t'-1)x+(d_{t'}-1)j\leq k-1 & \text{for all $t'\in[t]$}\\
&&& x,j \geq 0
\end{array}
\]
and note that as $xj-\lfloor x\rfloor\lfloor j\rfloor\leq x+j=O(k)$, we have $f^{(4)}_k(H)-O(k)\leq f_k^{(3)}(H)\leq f^{(4)}_k(H)$. Finally, by scaling $x$ and $j$ by $(k-1)$, we define
\[
\begin{array}{cccll}
c_H:=\frac{1}{(k-1)^2}
f_{k}^{(4)}(H) & = & \max & jx &\\
&& \text{s.t.} & (t'-1)x+(d_{t'}-1)j\leq 1 & \text{for all $t'\in[t]$}\\
&&& x,j \geq 0
\end{array}
\]
which is independent of $k$ and depends only on the graph $H$.
As $\rk_H^*(k)=f_k^{(4)}(H)\pm O(k)$, we finally attain $\rk_H^*(k)=\bigl(c_H+o(1)\bigr)k^2$.

Now, although the program for $c_H$ is not linear, it is clearly solvable in polynomial time. Furthermore, notice that for $(x,j)=\bigl({1\over 2(t-1)},{1\over 2(d_1-1)}\bigr)$, we have
\[
(t'-1)x+(d_{t'}-1)j\leq {1\over 2}+{1\over 2}=1,
\]
for all $t'\in[t]$, so $c_H\geq{1\over 4(t-1)(d_1-1)}$. Additionally, only considering the constraints $(1-1)x+(d_1-1)j\leq 1$ and $(t-1)x+(d_t-1)j\leq 1$, we find that $x\leq{1\over t-1}$ and $j\leq{1\over d_1-1}$, so $c_H\leq{1\over (t-1)(d_1-1)}$.
\end{proof}

\section{Conclusion and Further Directions}\label{sec:conclusion}

In our study of the extremal function $\rk_\fh(k)$, the largest open question is whether or not $\rk_\fh(k)=\rk_\fh^*(k)$ when $\fh$ is a family of simple, $r$-uniform graphs (see Conjecture~\ref{conj:multi}); also very natural is the question of the behavior of $\rk_{C_4}(k)$ (also discussed in the Introduction). We've noted that $\Omega(k^{4/3})\leq\rk_{C_4}(k)\leq\rk_{C_4}^*(k)\leq O(k^{3/2})$.
\begin{question}
	Is $\rk_{C_4}(k),\rk_{C_4}^*(k)=\Theta(k^{4/3})$?
\end{question}

Several further questions follow naturally from our line of inquiry.  For example:
\begin{question}
What are the exact asymptotics of $\rk_{P_t}(k)$?
\end{question}
 We note that, for a fixed $t$, $\rk_{P_t}(k)=\Theta(k^2)$, so we are interested in pinning down the constants here. More specifically, what are the extremal graphs for $\rk_{P_t}(k)$? Corollary~\ref{cor:connect} implies that there are extremal graphs for $\rk_{P_t}(k)$ with diameter at most $t$. Gy\'arf\'as, Rousseau and Schelp~\cite{GRS84} prove that if $n$ is sufficiently large compared to $t$, then 
 \[
 \ex(K_{n,n},P_t)=\begin{cases}
 {t-1\over 2}(2n-t+1) & \text{for $t$ odd};\\
 {t-2\over 2}(2n-t+2) & \text{for $t$ even.}
 \end{cases}
\]
This implies that for all $t\geq 5$ and $n$ sufficiently large, $\ex(K_{\sqrt{2}n},P_t)<\ex(K_{n,n},P_t)$, despite having the same number of edges, so most likely, the extremal graphs for $\rk_{P_t}(k)$ look more similar to cliques, as we showed was the case with $P_3$. However, $\ex(K_{\sqrt{2}n},P_4)\approx{3\over\sqrt{2}}n>2n\approx\ex(K_{n,n},P_4)$, so it may very likely be the case that the extremal graphs for $\rk_{P_4}(k)$ are bipartite. As there is this discrepancy, it would be very interesting to just determine the extremal graphs for $\rk_{P_4}(k)$ and why $P_4$ may behave differently from $P_t$ for all other $t$.

Next, in regard to the dumbbell $\cal{D}$ (see Theorem~\ref{thm:necklace}), we found that there is a multigraph $G$ on $(\phi-o(1))k$ edges with $\ex(G,\cal{D})<k$, but whenever $G'$ is a multigraph with $\ex(G',\cal{D})<k$ where each vertex has at most one loop, then $e(G)\leq{3\over 2}k<\phi\cdot k$. As such, it seems natural to ask about how $\rk_H^*(k)$ changes if $H$ is a non-uniform graph and the edges of different uniformities are weighted differently to reflect the fact that there are more possible edges of uniformity $2$ in the host graph than there are of uniformity $1$: one could more generally define $\ex(G,H):=\max\{\sum_{e\in E(F)}w(|e|): F \subset G, F ~ H\text{-free}\}$, where $w$ is an arbitrary weighting of the uniformities.

\begin{question}
How do $\rk_H(k)$ and $\rk_H^*(k)$ vary with the weight $w$ for non-uniform graphs?
\end{question}

In fact, since the main obstacle to forcing non-uniform graphs appears to be the edges having irreconcilable ``types,'' which suggests asking equivalent questions in the uniform case by artificially enforcing  distinct edge-types on graphs that are already uniform.

\begin{question}
Suppose $H$ is a graph consisting of both red and blue edges. How many edges can a red-blue colored graph $G$ have such that any subgraph with $k/2$ blue edges and $k/2$ red edges contains a colored copy of $H$?
\end{question}

We add the stipulation of having the same number of red and blue edges in light of how the Tur\'an question for arbitrary 2-colored cliques was already considered by Diwan and Mubayi~\cite{diwan2007turan}. One could request more generally that any $k$-edge subgraph of $G$ contains $H$, as we had for $H=\cal{D}$.

Finally, recall that we originally defined $\rk_\fh(k)$ by deciding that a host graph being ``best at forcing'' meant optimizing specifically its edge count, but one could just as easily ask this for any other monotone graph parameter $P$. That is, we could study $\rk_{P,\fh}(k):=\sup\{P(G):\ex(G,\fh)<k\}$. One particularly interesting example may be when $P=\chi$, the chromatic number. In this case, $\rk_{\chi,K_{1,t}}(k)$ and $\rk_{\chi,tK_2}(k)$ are not trivial.

\begin{question}
If $\fh$ is a family of (multi)(hyper)graphs, what is $\rk_{\chi,\fh}(k)$? 
\end{question}

To this end, we quickly note that as any 2-uniform graph $G$ has $e(G)\geq{\chi(G)\choose 2}$, Theorem~\ref{thm:esrank} implies that if $\fh$ is a family of simple, 2-uniform graphs with $\rho(\fh)=\rho\geq 3$, then $\rk_{\chi,\fh}(k)=\sqrt{\bigl(2+{2\over\rho-2}+o(1)\bigr)k}$; so again, it is most interesting to focus on families of bipartite graphs.

\section*{Acknowledgements}
We would like to thank Wesley Pegden for suggesting this problem and for helpful discussion. We would also like to thank Linhua Feng for spotting the relevance of ~\cite{AHS72}. Finally, we are very grateful to the anonymous referees whose comments have significantly improved this paper.

\end{document}